\documentclass[preprint,12pt]{elsarticle}

\usepackage{amssymb}
\usepackage{epstopdf}
\usepackage{amsfonts}
\usepackage{amsthm}
\usepackage{amsmath}
\usepackage{mathtools,xparse}
\usepackage[utf8]{inputenc}
\usepackage{url}
\usepackage{matlab-prettifier}
\usepackage{marvosym}
\usepackage{bm}
\usepackage[font=small,labelfont=bf]{caption}
\usepackage{float}
\usepackage{tikz}
\usetikzlibrary{quotes,angles,positioning}
\usepackage[ruled]{algorithm2e}
\usepackage{comment}
\usepackage{mathrsfs}
\usepackage{stmaryrd}
\usepackage[a4paper, left=4cm, right=2.5cm, top=1.5cm]{geometry}
\usepackage{pgfplots}
\usepackage{caption}
\usepackage{subcaption}
\usepackage{tikzsymbols}

\usetikzlibrary{matrix}
\usepgfplotslibrary{groupplots}
\pgfplotsset{compat=newest}
\usepackage{bigints}
\pgfplotsset{
    every axis plot/.append style = {font = \small}
  }

\newtheorem{problem}{Problem}[section]
\newtheorem{remark}{Remark}
\newtheorem{assumption}{Assumption}

\newcommand{\dof}{\text{dof}}
\newcommand{\Ndof}{N^{\dof}}

\newcommand{\materialDerivative}{\mathcal{D}_t}
\newcommand{\velocity}{\mathbf{v}}
\newcommand{\omegat}{\Omega^t}

\journal{Computers \& Mathematics with Applications}

\begin{document}

\begin{frontmatter}



\title{A Velocity-Based Moving Mesh Virtual Element Method}


\author[inst1]{H. Wells \corref{cor1}}
\cortext[cor1]{Corresponding author}
\ead{harry.wells@nottingham.ac.uk}

\affiliation[inst1]{organization={School of Mathematical Sciences, University of Nottingham},
            addressline={University Park}, 
            city={Nottingham},
            postcode={NG7 2RD}, 
            country={United Kingdom}}

\author[inst1]{M.E. Hubbard}
\author[inst2]{A. Cangiani}

\affiliation[inst2]{organization={Mathematics Area, SISSA},
            addressline={Via Bonomea 265}, 
            city={Trieste, 34136},
            country={Italy}}

\begin{abstract}
    We present a velocity-based moving mesh virtual element method for the numerical solution of PDEs involving boundaries which are free to move. The virtual element method is used for computing both the mesh velocity and a conservative Arbitrary Lagrangian-Eulerian solution transfer on general polygonal meshes. The approach extends the linear finite element method to polygonal mesh structures, achieving the same degree of accuracy. In the context of moving meshes, a major advantage of the virtual element approach is the ease with which nodes can be inserted on mesh edges. Demonstrations of node insertion techniques are presented to show that moving polygonal meshes can be simply adapted for situations where a boundary encounters a solid object or another moving boundary, without reduction in degree of accuracy.
\end{abstract}

\begin{keyword}
Moving Mesh Method \sep Virtual Element Method \sep Arbitrary Lagrangian-Eulerian Schemes \sep Polygonal Meshes \sep Porous Medium Equation
\PACS 0000 \sep 1111
\MSC 0000 \sep 1111
\end{keyword}

\end{frontmatter}


\section{Introduction}\label{sec::Intro}

We propose a new numerical method for  the simulation of free boundary problems over general polygonal meshes based on combining a  moving mesh algorithm with the Virtual Element Method (VEM). The new approach is presented focusing on the solution of the free boundary problem for the porous medium equation as a model problem. To showcase its generality,  we also briefly discuss an extension to a nonlinear forth-order problem.

Moving mesh methods form part of a large class of adaptive mesh refinement techniques alongside $h$- and $p$-refinement strategies. The primary advantage of moving mesh methods is the ability to optimize mesh structures without requiring any change in the mesh connectivity, which can cause computational challenges, particularly when parallel implementations are considered. They also provide a natural framework for tracking physical features of a time-dependent PDE, such as blow-up problems \cite{blowup}, phase change modelling \cite{phasechange}, fluid-structure interaction problems \cite{Richter2017}, and more general time-dependent PDEs \cite{baines1994moving,miller1981moving,gelinas1981moving}.

Here we consider a \emph{velocity-based} moving mesh algorithm \cite{MMFEM,MMFEMscale,MovingPhaseChange,MMFEMBC,MMFEMmonitors,MovingReview} for the numerical solution of  free boundary problems. The method is closely related to the Geometric Conservation Law method (GCL) \cite{GCL} and also forms part of a larger family of adaptive moving mesh methods \cite{MovingMeshesBook}. It uses a Lagrangian formulation of the given PDE to solve directly for the mesh velocities which are integrated over time to evolve the mesh. The solution computed on any given time-step is then transferred to the following time-step using an Arbitrary Lagrangian Eulerian (ALE) scheme \cite{ALEreview} based on a weak distribution of a given monitor function.

We combine the velocity-based moving mesh algorithm with the virtual element method for spatial discretisation on general polygonal meshes~\cite{basicprinciples}. The VEM was first introduced  in~\cite{basicprinciples} for the conforming discretisation of linear elliptic problems; cf. also~\cite{Eqproj,hitchiker,cangiani2017conforming,VEMMATLAB} for extensions and implementation details.
It provides a flexible discretization framework for the design of \emph{compatible} general mesh methods that incorporate, at the discrete level, fundamental properties of the continuous problem at hand, such as topology, conservation, symmetry, and positivity.
To achieve compatibility, discrete virtual element spaces are typically defined on individual polygonal elements implicitly through local boundary value problems. These are understood through a set of degrees of freedom instead of explicit basis functions. The use of local polynomial projections of the virtual discrete functions permits the definition of virtual element weak formulations which are computable at a cost on par with that of standard Finite Element Methods (FEMs). Here we consider the lowest-order VEM, which can be seen as a generalisation of the standard linear FEM to general polygonal/polyhedral meshes. The benefit of this generalisation is that typical issues such as node tangling, contact with obstacles, and topological changes of the moving domain~\cite{MovingMeshesBook,MovingReview} can be dealt with fully by local changes in the mesh topology. To demonstrate this,  we present a simple node insertion/removal algorithm for problems involving collision between a moving boundary and fixed obstacles. Just like linear FEM, the lowest-order VEM can be set up using only vertex information and in particular without the use of computationally expensive quadrature on polygons.
As such, our approach can be interpreted as a generalisation of the linear finite element moving mesh method~\cite{MovingReview}.

Numerical results demonstrate that optimal orders of convergence are achieved when applying a virtual element discretization on a diverse range of polygonal mesh structures. 
 To simplify the presentation, the novel moving mesh VEM is developed for the two-dimensional porous medium equation (PME). The PME is a parabolic non-linear diffusion equation \cite{PME}. It is an ideal benchmark for a moving mesh method as it provides time-varying solutions with compact support and a moving boundary. However, the moving mesh VEM proposed here may be immediately generalised to other free boundary problems.  Generalisations of the velocity-based method to other PDEs can be found in the literature \cite{MMFEM,MovingReview} and the VEM has already been generalised to the solution of a wide range of PDEs \cite{ellipticVEM,parabolicVEM,hyperVEM}, including flow in porous media \cite{Vacca2018,Zhao2020,Wang2019}. Hence we expect our approach to be extendable to a range of problems. To hint to the generality of the approach, in the numerical examples section we demonstrate the application of the moving mesh VEM  to a fourth-order diffusion problem with a moving boundary.

In recent years, several works have been published which involve time-dependent and moving polygonal domains. These developments include: ALE maps between two meshes using virtual element and discontinuous Galerkin (dG) methods \cite{dGALEmapping}, a cut-based dG scheme for fluid-structure interaction problems \cite{FSIPolyDG}, an ALE scheme on time-dependent reconstructed Voronoi meshes \cite{ALEmovingVoronoi}, adaptive methods using elliptic reconstruction techniques on moving meshes \cite{NonHeirVEM}, polygonal mesh quality measures for an anisotropic MMPDE scheme \cite{MMPDEAnisotropic}, and a VEM for long term simulations of moving landforms \cite{SedimentaryFlowVEM}. The benefits of using polygonal discretizations for moving mesh methods include representing moving boundaries and interfaces with a minimal number of degrees of freedom and localized mesh refinement when a change in mesh connectivity is required, such as in contact problems. To showcase the potential of the moving mesh VEM in this respect, we  test the method for the numerical treatment of both obstacle and self-intersecting  moving interface problems. 

The layout of the remaining sections of this paper is as follows. In Section \ref{sec::MovingMeshMethod} we introduce the moving mesh algorithm, including the weak formulations necessary for the method. The components of the VEM discretisation are split into two sections. Section \ref{sec::VEMVelocity} reviews the fundamentals of the VEM and formulates a method for computing mesh velocities at a fixed point in time. Section \ref{sec::VEMSolution} presents a framework for moving virtual elements along with a VEM for a conservative ALE update scheme. An overview of the algorithm and implementation details is given in Section \ref{sec::Implementation}. Node insertion/removal algorithms are discussed in Section \ref{sec::Contact}. A set of numerical experiments is presented in Section \ref{sec::Tests} for the PME and a fourth-order diffusion problem before drawing some conclusions in Section~\ref{sec::conclusions}.

\section{The Moving Mesh Framework}\label{sec::MovingMeshMethod}

In this section we present the key components of the velocity-based moving mesh framework for a general time-dependent PDE of the form
\begin{equation} \label{eq::PDE}
  \frac{\partial \rho}{\partial t} = \mathcal{L} \rho,
\end{equation}
with $t\in (0,T]$, $T\in\mathbb{R}^+$, indicating time and $\mathcal{L}$ a generic spatial differential operator in $\mathbb{R}^d$. The choice of spatial operator will be made specific later on when we introduce the porous medium equation as the model problem. Further details on a finite element discretisation of this moving mesh method for general parabolic problems can be found in \cite{MovingReview} and the references therein.

The time-dependent support of the evolving solution to Equation \eqref{eq::PDE} will be denoted by $\Omega^t \subseteq \mathbb{R}^d$, with boundary $\partial \Omega^t$ which is allowed to move.  
The temporal superscript $t$ will be swapped for $n$ when considering discrete time levels $t^n$ or, when no confusion arises, the superscript notation will be dropped entirely for ease of reading. 

A moving coordinate system is defined by $\mathbf{x} = \mathbf{x}(\boldsymbol{\xi},t)$ where $\mathbf{x}(\boldsymbol{\xi},0) = \boldsymbol{\xi}$ is the coordinate system on the initial  domain $\Omega^0$, which is used as the reference domain for simplicity. 
It is assumed that, for all $t\in [0,T]$, the map $\mathbf{x}(\cdot,t)$ is Lipschitz with Lipschitz inverse. 
The evolution of $\Omega^t$ is determined by the velocity field
\begin{equation} \label{eq::VelocityODE}
    \velocity = \frac{\partial \mathbf{x}}{\partial t}.
\end{equation}
Following~\cite{Bonito2013}, we shall also refer to the space-time domain
\[
\mathcal{Q}_T=\{(\mathbf{x},t)\, : \, t\in [0,T], \mathbf{x} = \mathbf{x}(\boldsymbol{\xi},t), \boldsymbol{\xi}\in \Omega^0 
\}.
\]

\subsection{The Velocity Field}

The movement of the mesh is derived by specifying the time evolution of the distribution of the spatial integral of some solution-dependent monitor function, $\mathbb{M}(\rho)$. Specifically, the velocity field $\mathbf{v}$
is determined by requiring that the initial distribution of the monitor function $\mathbb{M}(\rho)$ is conserved as time progresses. 
Namely, we look for $\mathbf{v} = \frac{\partial \mathbf{x}}{\partial t}$ such that, for all $w \in L^2(\mathcal{Q}_T)$, the coordinates $\mathbf{x}(\boldsymbol{\xi},t)$ satisfy
\begin{equation} \label{eq::Distribution}
 \frac{\int_{\omegat} w(\mathbf{x},t) \, \mathbb{M}(\rho(\mathbf{x},t))\ d \mathbf{x}}
      {\int_{\omegat} \mathbb{M}(\rho(\mathbf{x},t))\ d \mathbf{x}} =
 \frac{\int_{\Omega^0} w(\boldsymbol{\xi},0) \, \mathbb{M}(\rho(\boldsymbol{\xi},0))\ d \boldsymbol{\xi}}
      {\int_{\Omega^0} \mathbb{M}(\rho(\boldsymbol{\xi},0))\ d \boldsymbol{\xi}}
                   \qquad\quad \forall \, t\in [0,T].
\end{equation}
The derivation of the moving mesh method is based on application of the Reynolds transport theorem and on the assumption that the material derivatives of the weighting functions $w(\mathbf{x},t)$ are zero with respect to the velocity field $\velocity$ (see Equation \eqref{eq::VelocityODE}). Namely,
\begin{equation} \label{eq::MaterialDerivative}
    \materialDerivative w = \frac{\partial w}{\partial t} + \velocity \cdot \nabla w = 0 .
\end{equation}
This is a common assumption made in finite element approaches to moving mesh algorithms \cite{dGALEmapping,Bonito2013,Bonito2013apriori} and is equivalent to assuming that $w(\mathbf{x(\boldsymbol{\xi},t)},t) = w(\boldsymbol{\xi},0)$ in Equation \eqref{eq::Distribution}. 

Equidistribution-based mesh movement algorithms typically attempt to reduce the global approximation error without changing the number of degrees of freedom by choosing a finite set of weighting functions $w_i(\mathbf{x},t)$, so that each one has a direct association with a mesh node or element. The monitor function $\mathbb{M}(\rho)$ is then selected to act as a local error indicator, and the mesh is adjusted in an attempt to equidistribute the values of the weighted monitor integrals
\begin{equation} \label{eq::WeightedMonitorIntegral}
 \mu^t(w_i) = \int_{\omegat} w_i \, \mathbb{M}(\rho) \ d \mathbf{x},
\end{equation}
and hence to equidistribute the local error across the mesh. Such an approach could be followed within our framework \cite{MMFEMmonitors,MovingReview}. However, in this paper we adopt the following approach which is more akin to Lagrangian mesh movement algorithms, which attempt to move the mesh with the `flow' velocity.

Choosing $\mathbb{M}(\rho) = \rho$ in equation (\ref{eq::WeightedMonitorIntegral}) naturally leads to a weak approximation of the Lagrangian `flow' velocity of the PDE when $\mathcal{L} \rho = \nabla \cdot \mathbf{f}$ for some flux $\mathbf{f}$ in Equation \eqref{eq::PDE}. This has two benefits: (a) it allows us to predict the movement of free boundaries; (b) it reduces interpolation error of the computed mesh and solution between time-steps because the mesh (and hence the solution) is transported with the velocity field inherent to the PDE. 

In many PDEs (including the PME), $\rho$ represents a density, \emph{i.e.}\ its integral in space is a mass, so we will refer to $\mathbb{M}(\rho) = \rho$ as the mass monitor. For clarity of presentation, we will assume use of this mass monitor from now on. Hence, each weight function $w(\mathbf{x},t)$ is assigned its own `mass':
\begin{equation} \label{eq::WeightedMasses}
  \mu^t(w) = \int_{\omegat} w \rho\ d \mathbf{x}
           = c(w) \theta^t \qquad\quad \mbox{where} \qquad\quad
  \theta^t = \int_{\omegat} \rho\ d \mathbf{x} .
\end{equation}
Equation \eqref{eq::Distribution} then provides $c(w) = \mu^0(w) / \theta^0$ from the initial conditions of the PDE, and we assume that these distribution coefficients remain constant in time. The evolution of $\rho$ is governed by Equation \eqref{eq::PDE}, so Equation \eqref{eq::WeightedMasses} provides us with a way to prescribe the evolution of the coordinate system in a way which retains the initial `mass' distribution.
\

\begin{remark}
The original moving mesh finite element method \cite{MMFEM,MovingReview} chose the weight functions to be the standard linear Lagrange finite element test functions on meshes of simplices. This associates a fixed proportion of the total mass of the system with each node of the mesh, so the values of $\mu^t(w)$ depend not only on $\rho$, but also on the mesh node positions. In this situation, Equation \eqref{eq::WeightedMasses} provides a way to compute mesh node velocities using standard finite element techniques, in a way which is consistent with local mass conservation when this is a property of the underlying PDE. Our results to follow demonstrate that the same principle can be applied on polygonal meshes within a virtual element framework.
\end{remark}

\

Let $t\in [0,T]$ and consider  all $w\in H^1(\mathcal{Q}_T)$ such that $\materialDerivative w(t) =0$, cf. Equation \eqref{eq::MaterialDerivative}.  Differentiating the first equality in~\eqref{eq::WeightedMasses} with respect to time and applying the Reynolds transport theorem~\cite{CFD} leads to
\begin{equation}\label{eq::MudotMH}
    \dot{\mu}^t(w) = \int_{\omegat} w \left\{ \frac{\partial \rho}{\partial t}
        + \nabla \cdot \left( \rho \velocity \right)\right\}\ d \mathbf{x} .
\end{equation}
Noting  that~\eqref{eq::MudotMH} does not fully determine the velocity, we further require  that  each weighting function retains a fixed proportion of $\theta^t$ as the mesh evolves. Hence we impose  the constraints
\begin{equation} \label{eq::MovingConservationLaw}
   \dot{\mu}^t(w) = c(w) \dot{\theta}^t.
\end{equation}
Inserting~\eqref{eq::MovingConservationLaw} in~\eqref{eq::MudotMH}, using~\eqref{eq::PDE}, and integrating by parts  results in
\begin{eqnarray}
  c(w) \dot{\theta}^t\ +\ \int_{\omegat} \rho \nabla w \cdot \mathbf{v}\ d\mathbf{x}
     &=& \int_{\omegat} w \mathcal{L} \rho\ d\mathbf{x}\ 
        +\ \int_{\partial \omegat} w \rho \mathbf{v} \cdot \mathbf{n}\ ds, \label{eq::GeneralVelocity1}
\end{eqnarray}
with $\mathbf{n}$ the outward unit normal on $\partial \omegat$.
Similarly, a direct application of the Reynolds transport theorem to the second equality in~\eqref{eq::WeightedMasses}, yields
\begin{eqnarray}
  \dot{\theta}^t &=& \int_{\omegat} \mathcal{L} \rho\ d\mathbf{x}\
        +\ \int_{\partial \omegat} \rho \mathbf{v} \cdot \mathbf{n}\ ds.\label{eq::GeneralVelocity2}
\end{eqnarray}

In fact, $\dot{\theta}^t$ is typically known explicitly because $\rho \mathbf{v} \cdot \mathbf{n}$ is known on the whole of $\partial \omegat$, and the boundary conditions provided with the PDE will enable evaluation of the integral of $\mathcal{L} \rho$. Furthermore, for mass-conservative problems, $\dot{\theta}^t = 0$. Note also that assuming that $c(w)$ remains constant preserves the initial distribution of the mass (or a more general monitor integral) so a monitor which is initially equidistributed between a set of weighting functions should remain equidistributed as the coordinate system evolves with this velocity field.

Equations~\eqref{eq::GeneralVelocity1} and~\eqref{eq::GeneralVelocity2} are used to compute an instantaneous  velocity $\velocity$ which is consistent with conserving the proportion of mass associated with each weighting function. This provides a form of local mass conservation when $\dot{\theta}^t = 0$. However, it still does not uniquely define $\mathbf{v}$ in multiple space dimensions. 
To overcome this issue, the velocity field is written in terms of its Helmholtz decomposition
\begin{equation}\label{eq::VelocityDecomp}
    \mathbf{v} = \mathbf{q} + \nabla \phi,
\end{equation}
where $\phi$ is a scalar potential and $\mathbf{q}$ must be specified. This constraint is equivalent to imposing the curl of the velocity field because $\nabla \times \mathbf{q} = \nabla \times \mathbf{v}$. Moreover, for simplicity, we may further assume that $\mathbf{q} = \mathbf{0}$. This is the natural choice for the porous medium equation, which is derived under the assumption of a curl-free flow velocity field. An example of the method applied with a rotational velocity field is presented in \cite{MovingReview}. The problem for determining the velocity potential is therefore: find $\phi\in H^1(\omegat)$ such that
\begin{equation} \label{eq::VelocityPotential}
  \int_{\omegat} \rho \nabla w \cdot \nabla \phi \ d\mathbf{x}
     = \int_{\omegat} w \mathcal{L} \rho\ d\mathbf{x} \ 
        +\ \int_{\partial \omegat} w \rho \mathbf{v} \cdot \mathbf{n}\ ds \
        -\ c(w) \dot{\theta}^t \qquad\forall w\in H^1(\omegat),
\end{equation}
where $\phi = 0$ is specified at one point in $\omegat$ to ensure uniqueness.

The velocity field is finally obtained as the solution of equation~\eqref{eq::VelocityDecomp} written in weak form with  $\mathbf{q} = \mathbf{0}$ and $\phi$ given by~\eqref{eq::VelocityPotential}. That is, we find $\velocity\in [H^1(\omegat)]^d$ such that
\begin{equation} \label{eq::VelocityRecovery}
 \int_{\omegat} z {\velocity} \ d\mathbf{x} = \int_{\omegat} z \nabla \phi \ d \mathbf{x} \qquad\quad \forall \, z \in H^1(\omegat),
\end{equation}
with $\velocity \cdot \mathbf{n}$ imposed on any part of the boundary where it is known. For instance, in case of contact with an obstacle, we would impose  $\velocity \cdot \mathbf{n}=0$ on the contact boundary. Note that~\eqref{eq::VelocityRecovery} is just the component-wise $L^2$-projection.

\begin{remark}
The velocity field in the interior of $\omegat$ could be discarded at this stage, and replaced by one derived using a different approach which is constrained by the boundary velocity derived above. For example, interior mesh nodes could be moved using our approach with a different monitor function or a Laplacian smoothing could be applied \cite{MMFEMmonitors,Richter2017}.
\end{remark}

\subsection{Solution Update}\label{sec:solupdate}

The velocity field ${\velocity}$ derived using Equations \eqref{eq::VelocityPotential} and \eqref{eq::VelocityRecovery} can now be used in the update of the solution. Integration of Equation \eqref{eq::MudotMH} by parts and substitution of Equation \eqref{eq::PDE} results in
\begin{equation} \label{eq::ALEUpdate}
  \dot{\mu}^t(w) = \int_{\omegat} w \mathcal{L} \rho\ d\mathbf{x}\ 
        -\ \int_{\omegat} \rho \nabla w \cdot {\mathbf{v}}\ d\mathbf{x}\
        +\ \int_{\partial \omegat} w \rho \mathbf{v} \cdot \mathbf{n}\ ds
           \qquad\quad \forall \, w \in H^1(\omegat) .
\end{equation}
This is a standard conservative ALE update. Along with Equation \eqref{eq::VelocityODE}, it gives a system of ODEs governing the evolution of the coordinate system and the distribution of the mass monitor which can be approximated using standard solvers such as explicit Runge-Kutta methods~\cite{MovingReview,MMFEMBC}. With these at hand, the solution can be recovered by solving the problem: find $\rho\in H^1(\omegat)$ such that
\begin{equation} \label{eq::SolutionRecovery}
  \int_{\omegat} w \rho \ d \mathbf{x} = \mu^t(w) \qquad\quad \forall \, w \in H^1(\omegat) .
\end{equation}
 In the special case where $\dot{\theta} = 0$ and $c(w)$ is assumed constant in time, Equation \eqref{eq::ALEUpdate} is redundant because $\dot{\mu}^t(w) \approx 0$ by design. In fact, \emph{without} the velocity recovery from the potential through~\eqref{eq::VelocityPotential} and~\eqref{eq::VelocityRecovery}, we would have $\dot{\mu}^t(w) \equiv 0$. However, the practical steps of the recovery procedure may introduce small perturbations. Alternatively, one may use the knowledge that $\dot{\mu}^t(w) \equiv 0$ and hence resort directly to the initial values ${\mu}^0(w)$, avoiding the need to calculate the ALE update \eqref{eq::ALEUpdate} altogether: we refer to this as direct recovery. However, direct recovery can only be used with the specific choice of the mass monitor and with mass-conservative PDEs. In other situations the interior velocity field will not correspond to $\dot{\mu}^t(w) = 0$ and the ALE update is essential.

\begin{remark}
An alternative, non-conservative, ALE formulation can be derived which would give an equation for $\dot{\rho}$ instead of $\dot{\mu}$ (see \cite{MMFEMmonitors}). This derivation does not require that the weight functions evolve with zero material derivative, but conservation of mass is lost as a result.
\end{remark}

\subsection{Porous Medium Equation}
On an open, bounded domain $\tilde{\Omega} \subset \mathbb{R}^d$, we consider the following initial-boundary value problem for the  Porous Medium Equation (PME): find $\rho:\tilde{\Omega} \times\mathbb{R}^+\rightarrow\mathbb{R}$ such that
\begin{align}
\label{eq:PME1}
    \frac{\partial \rho}{\partial t}(\mathbf{x},t) &= \Delta \Phi(\rho(\mathbf{x},t)) \qquad && (\mathbf{x},t)\in \tilde{\Omega}\times (0,\infty],\\
    \label{eq:PME2}
     \rho(\mathbf{x},t) &= 0 \qquad &&(\mathbf{x},t)\in\partial \tilde{\Omega}\times (0,\infty], \\
     \label{eq:PME3}
     \rho(\mathbf{x},0) &= \rho^0(\mathbf{x})\qquad &&\mathbf{x} \in \tilde{\Omega},
\end{align}
where $\Phi= \rho^{m+1}/(m+1)$, for some $m > 0$, and $\rho^0\ge 0$ having compact support in $\tilde{\Omega}$. The PME belongs to the broader class of Generalized Porous Medium Equations (GPME), also known as filtration equations, obtained with  $\Phi:\mathbb{R}^+\rightarrow\mathbb{R}^+$ any increasing function. 
The mathematical analysis of the GPME is well developed; see the monograph~\cite{PME} and the references therein. In particular, the notion of appropriate weak solutions is discussed in~\cite{PME} where it is shown that, for
non-negative $\rho^0 \in L^1(\tilde{\Omega})$ and for $\Psi(\rho^0) \in L^1(\tilde{\Omega})$, where $\Psi$ is the anti-derivative of $\Phi$,  if $\Psi(\rho)>0$ for $\rho>0$,  there exists a unique non-negative weak solution to the GPME globally in time.

 The PME models a number of  physical processes such as fluid flow, heat transfer, and diffusion. It exhibits several interesting properties, including the existence of a family of radially symmetric similarity solutions \cite{MathBioBook}, which are used to test the numerical method in Section \ref{sec::Tests}. Other properties of the PME are discussed in \cite{PME,MovingReview,MovingHessian}.

It is a fundamental example of a degenerate parabolic equation, stemming from the condition that $\Psi$ is non-negative, rather than simply positive.

Solutions that exhibit an evolving compact support are ideal for the class of moving mesh methods considered herein because considering as unknown the support of the solution leads to a moving boundary problem that can be simulated over time without having to discretize the entire geometry of $\tilde{\Omega}$. 

Introducing the time-dependent support of $\rho$ as $\omegat$, we define the time-dependent coordinate system $\mathbf{x}(\boldsymbol{\xi},t)$ with a velocity field $\dot{\mathbf{x}} = \mathbf{v}(\mathbf{x},t)$ that corresponds, at $\partial\Omega^t$, with the movement of the free boundary of $\omegat$ for all $t \in [0,\infty]$.
Additionally, we consider the free boundary problem for the PME in which part of the boundary may be obstructed by a fixed object. To this end, the boundary is divided into a moving part $\partial \omegat_M$ and a fixed part $\partial \omegat_F$, such that $\partial \omegat = \partial \omegat_M \cup \partial \omegat_F$ and $\partial \omegat_M \cap \partial \omegat_F = \emptyset$. Thus, we arrive to the following classical free boundary problem for the PME, whose smooth solutions are weak solutions of~\eqref{eq:PME1}-\eqref{eq:PME3}, cf.~\cite{PME}.
\begin{problem}[The Porous Medium Equation (PME)] \label{prob::PME}~\\
Let $T>0$ and $m > 0$. Find $\rho = \rho(\mathbf{x},t)$ such that
$\rho(\mathbf{x},0) = \rho^0(\mathbf{x})$ for $\mathbf{x} \in \Omega^0$ and,
for all $t \in (0,T]$,
\begin{align*}
    \frac{\partial \rho}{\partial t} &= \nabla \cdot (\rho^m \nabla \rho) &&\mathbf{x} \in \omegat, \\
    \rho &= 0 &&\mathbf{x} \in \partial \omegat_M, \\
    \rho^m \nabla \rho \cdot \mathbf{n} &= 0 &&\mathbf{x} \in \partial \omegat_F.
\end{align*} 
Here, $\mathbf{n}$ is the outward pointing unit normal to the boundary $\partial \omegat$. Note that an additional (kinematic) boundary condition,
\begin{equation*}
    \rho \mathbf{v} \cdot \mathbf{n} = - \rho^m \nabla \rho \cdot \mathbf{n} \qquad\qquad\qquad \mathbf{x} \in \partial \omegat_M,
\end{equation*}
which imposes zero flux through the moving boundary, is required to determine the boundary velocity $\mathbf{v}$. On $\partial \Omega_F$, the boundary velocity is specified and $\rho \mathbf{v} \cdot \mathbf{n} = 0$.
\end{problem}

Given a transformed coordinate system $\mathbf{x}$ and a distribution $\mu$ of the variable $\rho$, our method for solving the PME Problem~\ref{prob::PME} requires the solution of the following problems successively.

\begin{problem}[Solution Reconstruction -- Equation \eqref{eq::SolutionRecovery}] \label{prob::SolutionRecon}~\\
Given $\mu^t(w) \in \mathbb{R}$, find $\rho \in H^1(\omegat)$ such that
\begin{equation}\label{eq::ContinuousALEUpdate}
    m(\rho,w) = \mu^t(w) \qquad\quad \forall \, w \in H^1(\omegat) ,
\end{equation}
where
\begin{equation}
    m(\rho,w) = \int_{\omegat} w \rho \ d\mathbf{x}.
\end{equation}
\end{problem}

\begin{problem}[Velocity Potential -- Equation \eqref{eq::VelocityPotential}] \label{prob::ContinousPotential}~\\
Given $\rho \in H^1(\omegat)$, find $\phi \in H^1(\omegat) \slash \mathbb{P}_0(\Omega^t)$, where $\mathbb{P}_0(\Omega^t)$ is the space of constants, such that
\begin{equation} \label{eq::Potential}
    {A}(w,\phi) = d(w) \qquad\quad \forall \, w \in H^1(\omegat)\slash \mathbb{P}_0(\Omega^t),
\end{equation}
where
\begin{align}
    &{A}(\phi,w) = \int_{\omegat} \rho \nabla \phi \cdot  \nabla w\ d \mathbf{x}, \label{eq::A_w}\\
    &d(w) = -\int_{\omegat} \rho^m \nabla \rho \cdot \nabla w\ d\mathbf{x}.\label{eq::d_w}
\end{align}
\end{problem}
\noindent
This has been derived by substituting $\mathcal{L} \rho = \nabla \cdot (\rho^m \nabla \rho)$ in \eqref{eq::VelocityPotential} and applying the boundary conditions from Problem \ref{prob::PME}. To reconstruct the velocity field we introduce the follwoing space
\begin{equation}
\left[ H^1_{\partial \omegat_F}(\omegat) \right]^d = \left\{ \mathbf{w} \in \left[ H^1(\omegat) \right]^d\ :\ \mathbf{w} \cdot \mathbf{n}|_{\partial \omegat_F} = 0 \right\},
\end{equation}
in which we have the following reconstruction problem.

\begin{problem}[Velocity Reconstruction -- Equation \eqref{eq::VelocityRecovery}] \label{prob::VelocityReconstruction}~\\
Given $\phi \in H^1(\omegat)$, find ${\velocity} \in \left[ H^1_{\partial \omegat_F}(\omegat) \right]^d$ such that
\begin{equation}\label{eq::VelocityReconstruction}
    M({\mathbf{v}},z) = b(z) \qquad\quad \forall \, z \in H^1(\omegat),
\end{equation}

where
\begin{align}
    M({\mathbf{v}},z) &= \int_{\omegat} z {\velocity} \ d\mathbf{x}, \\
    b(z) &= \int_{\omegat} z \nabla \phi \ d\mathbf{x}.
\end{align}
\end{problem}
\noindent
To compensate for the fact that the velocity is reconstructed from the potential (see the comments in Section~\ref{sec:solupdate}), an \emph{a posteriorily} computed ALE update may be used~\cite{ALEreview}. 
Specifically, substituting $\mathcal{L} \rho = \nabla \cdot (\rho^m \nabla \rho)$ and the PME boundary conditions into~\eqref{eq::ALEUpdate} we arrive to the following.

\begin{problem}[ALE Update -- Equation \eqref{eq::ALEUpdate}] \label{prob::ALEUpdate}~\\
Given $\rho \in H^1(\omegat)$ and ${\velocity} \in \left[ H^1(\omegat) \right]^d$, fix
\begin{equation}\label{eq::ALEUpdateContinuous}
    \dot{\mu}^t(w) = \int_{\omegat} - \rho \nabla w \cdot \left\{ \rho^{m-1} \nabla \rho + {\velocity} \right\} \ d\mathbf{x} \qquad\quad \forall \, w \in H^1(\omegat).
\end{equation}
\end{problem}
\noindent
Once again we note that, when $\dot{\mu}^t(w) = 0$ is assumed, Equation \eqref{eq::ALEUpdateContinuous} is redundant and Equation \eqref{eq::Potential} is equivalent to a weak form of the PDE
\begin{equation}
  \nabla \cdot ( \rho \velocity ) = - \nabla \cdot ( \rho^m \nabla \rho ) ,
\end{equation}
where $\nabla \times \velocity = \mathbf{0}$ has been assumed. In other words, the formulation recovers the Lagrangian flow field associated with the porous medium equation. 

Equations \eqref{eq::ALEUpdateContinuous} and \eqref{eq::VelocityReconstruction} can now be combined to give a system of ODEs which determine the evolution of the coordinate system and the monitor distribution. These can be discretized in time using an appropriate solver for first-order ODEs. Specifically, the forward Euler method is used to perform the tests presented in Section~\ref{sec::Tests}. Higher-order  explicit time-stepping schemes have been presented in \cite{MMFEMBC,MovingReview}, but they are not required here.

We now proceed to describe how Problems \ref{prob::SolutionRecon}--\ref{prob::ALEUpdate} can be approximated in space on two-dimensional polygonal meshes using the virtual element method.

\section{Virtual Element Method for the Velocity}\label{sec::VEMVelocity}
In this section we introduce the VEM and apply it for the discretisation of Problems ~\ref{prob::ContinousPotential} and~\ref{prob::VelocityReconstruction} to determine the domain velocity field.
We follow the framework introduced in~\cite{basicprinciples,Eqproj,cangiani2017conforming} for the discretisation of linear elliptic problems. For simplicity, we consider the lowest-order (linear) VEM in two dimensions. Hence, for the rest of the paper, it is assumed that $d=2$. 

In the following we let $0=t^0<t^1<\dots t^{N_t}=T$ denote a given sequence of discrete time levels at which the PME is discretised. 

\subsection{A Moving Polygonal Mesh}\label{sec::MeshAssumptions}
The polygonal  representation of the moving domain $\omegat$ at the discrete time level $t^n$ is denoted by $\Omega^n_h$ and the corresponding polygonal mesh partitioning $\Omega^n_h$ is denoted by $\mathcal{T}_h^n$. Polygonal elements and edges within the mesh are denoted by $E$ and $e$, respectively. For a given polygonal element $E$, $|E|$ denotes the area, $h_E$ the diameter, and $(x_c,y_c)$ the centroid.

The following conditions are sufficient to guarantee convergence of the VEM and are assumed to hold true at each time level $t^n$, $n=0,1,\dots, N_t$; see, for example, \cite{basicprinciples}.

\begin{assumption}[Mesh Partitioning]\label{assumption::MeshPartitioning}
 The mesh $\mathcal{T}_h^n$ 
provides a partition of $\Omega^n_h$ into non-overlapping simple polygons\footnote{A polygon is \emph{simple} if its boundary  forms a closed graph with no intersecting edges other than at each vertex which intersects exactly two edges.
}. 
\end{assumption}

\begin{assumption}[Shape Regularity]\label{assumption::ShapeRegular}
Every $E \in \mathcal{T}_h^n$ is a star-shaped domain or a finite union of star shaped domains with respect to a ball of radius greater than $\gamma h_E$ for some uniform $\gamma > 0$. Additionally,  for all edges $e \in \partial E$, the length of $e$ is greater than $\delta h_E$ for some  uniform $\delta > 0$.
\end{assumption}

\begin{remark}
The shape regularity assumption on the edges required by Assumption~\ref{assumption::ShapeRegular} can be entirely removed following, for example, \cite{brenner2018virtual}. Small edges are  relevant, for instance, for the various types of contact problems considered below.
\end{remark}

\subsection{The $\Pi^{\nabla}$ and $\Pi^0$ Projections}
Given a subset $\omega \subset \mathbb{R}^d$, $d=1,2$, representing an element or an edge,  with  $\mathbb{P}_k(\omega)$ we denote the space of polynomials of degree $k$ on $\omega$. 
Furthermore, we use $(\cdot,\cdot)_\omega$ to denote the $L^2$ inner product over $\omega$.

Critical to the construction of a VEM is the judicious choice of discrete spaces and degrees of freedom that permit the computation of local polynomial projections  of virtual functions using only the degrees of freedom.

Specifically, the definition of the VEM space to be introduced below is based on an $H^1$-type projection operator $\Pi^{\nabla}$ onto $\mathbb{P}_1(E)$ given, for any $w\in  H^1(E)$, by
\begin{equation}\label{eq::H1Projection}
\left\{
\begin{array}{l}
    \left(\nabla \Pi^{\nabla} w, \nabla p \right)_{0,E} = \left( \nabla w, \nabla p \right)_{0,E}\ \ \ \forall p \in \mathbb{P}_1(E)  \\ 
    \left(w-\Pi^{\nabla}_1 w,1  \right)_{0,\partial E} = 0. 
\end{array}
\right.
\end{equation}
Further, the VEM requires the availability of the $L^2$-projection operator $\Pi$ onto  $\mathbb{P}_1(E)$ defined, for $w\in  L^2(E)$, by
\begin{equation}\label{eq::L2Projection}
        \left( \Pi w, p \right)_{0,E} = \left( w,p \right)_{0,E}\ \ \ \forall p \in \mathbb{P}_1(E),
\end{equation}
Additionally, the VEM requires the $L^2$-projection onto constants of $\nabla w$ which we denote by  $\mathbf{\Pi}_0$, thus
 \begin{equation}\label{eq::GradProjection}
   \left(\mathbf{\Pi}_0\nabla w, \mathbf{p} \right)_{0,E} = \left( \nabla w,\mathbf{p} \right)_{0,E}\ \ \ \forall \mathbf{p} \in [\mathbb{P}_0(E)]^2.
\end{equation}
The same operator applied component-wise to vector-valued functions will be denoted by  $\mathbf{\Pi}$. 

\subsection{Virtual Element Spaces} \label{sec::LocalVEMSpace}
Virtual element discrete functions are defined implicitly as solutions of  element-wise boundary value problems. These functions are only accessed  through their degrees of freedom with which specific projections can be computed; in this way,  \emph{no} evaluation of the (implicitly known) discrete functions is required. 
Here we consider the linear virtual element space of~\cite{Eqproj} which is a modification of the original space proposed in~\cite{basicprinciples} for which both  the $H^1$ and $L^2$ projections introduced in the previous section are computable from the degrees of freedom.

Given a polygonal element $E \in \mathcal{T}_h^n$,  we first define the elemental boundary space as
\begin{equation}\label{eq::BoundaryVEMSpace}
    \mathbb{B}(\partial E) = \{ w_h \in C^0(\partial E)\ :\ w_h|_{e} \in \mathbb{P}_1(e)\ \ \ \forall e \subset \partial E \}.
\end{equation}
Then the local virtual element space of~\cite{basicprinciples}  is defined as
\begin{equation}\label{eq::VEMSpace}
    W(E) = \{ w_h \in H^1(E)\ :\ w_h|_{\partial E} \in \mathbb{B}(\partial E),\ \Delta w_h|_E = 0 \}.
\end{equation}
Note that $\mathbb{P}_1(E) \subseteq W(E)$. In fact, it is immediate to check that the space $W(E)$ corresponds to the elemental linear finite element space whenever $E$ is a triangle, that is $W(E) \equiv \mathbb{P}_1(E)$. In this respect, the VEM can be seen as a generalisation of the linear FEM to polygonal meshes.
It is clear from the definition that the only degrees of freedom in the choice of a function $w_h\in W(E)$ are related to its values on the boundary and, given that $w_h$ is linear on each edge, ultimately they depend only on the nodal values; a rigorous proof that these constitute a unisolvent set of degrees of freedom  for $W(E)$ is presented in \cite{basicprinciples}.
Hence, in analogy to classical continuous linear finite element spaces, we identify the VEM functions by their nodal values.

Additionally, the nodal values allow for the computation of the local $\Pi^{\nabla}$ projection of any $w_h \in W(E)$. Indeed, in order to compute the first equation in~\eqref{eq::H1Projection} we need to be able to evaluate, for any $p \in \mathbb{P}_1(E)$, the right-hand side term:
\[
\left( \nabla w_h, \nabla p \right)_{0,E}=-\left( w_h, \Delta p \right)_{0,E}+\left( w_h, \mathbf{n}\cdot \nabla p \right)_{0,\partial E}=\sum_{e\in\partial E}\left( w_h, \mathbf{n}\cdot \nabla p \right)_{0,e},
\]
as $\Delta p=0$. This shows that such terms only depend on the value of $w_h$ on the boundary of $E$, which are known and are determined by the nodal values. The same is clearly true for the second equation in~\eqref{eq::H1Projection}, hence the  local $\Pi^{\nabla}$ projection is fully computable just by accessing the degrees of freedom. 
However, the $L^2$ projection is not computable for this space. Hence, following~\cite{Eqproj}, we first consider the  augmented space
\begin{equation}
    \tilde{V}(E) = \{ w_h \in H^1(E)\ :\  w_h|_{\partial E} \in \mathbb{B}(\partial E),\ \Delta w_h \in \mathbb{P}_1(E) \},
\end{equation}
from which a new \emph{enhanced} local virtual element space can be defined as
\begin{equation}\label{eq::NewVEMSpace}
    V(E) = \{ w_h \in \tilde{V}(E)\ :\ (w_h - \Pi^{\nabla} w_h, q)_{0,E} = 0\ \ \ \forall q \in \mathbb{P}_1(E) \}.
\end{equation}
Notice that the space $V(E)$ is characterised by the following four properties: (i) we still have that $\mathbb{P}_1(E) \subseteq V(E)$; (ii) the dimension of $V(E)$ is equal to that of $W(E)$ and we can still use the nodal values as degrees of freedom (see \cite{Eqproj});  (iii) the  local $\Pi^{\nabla}$ projection is still computable; 
(iv) if $w_h\in V(E)$ then $\Pi w_h=\Pi^{\nabla} w_h$ by construction, hence also the $L^2$-projection $\Pi$ is computable.
The proof that the gradient projection into constants $\mathbf{\Pi}_0\nabla w_h$ is also computable is similar.

The global virtual element space is then defined for a given time $t^n$ as
\begin{equation}\label{eq::GlobalVEMSpace}
    V_{h}^n = \{ w_h \in H^1(\Omega_h^n) : w_h|_E \in V(E)\ \ \ \forall E \in \mathcal{T}_h^n \}.
\end{equation}
The dimension of this space is equal to the number of nodes in the mesh, which we shall denote by $\Ndof$.

\begin{remark}
When required, homogeneous Dirichlet boundary conditions can be embedded in the virtual element space by fixing the relevant boundary nodes. Hence, given a $\Gamma \subseteq \partial \Omega_h^n$, we denote the constrained space by $V_{h,\Gamma}^n = \{ w_h \in V_h^n\ :\ w_h|_{\Gamma} = 0 \}$.
\end{remark}

\subsection{Discretisation of the velocity problems} \label{sec::VelocityForms}
Having defined the  virtual element space in Section \ref{sec::LocalVEMSpace}, we are ready to present the VEM for the solution of  the velocity problems \ref{prob::ContinousPotential} and \ref{prob::VelocityReconstruction}. The VEM is based on the construction of approximate weak forms which are \emph{computable} through the elemental projection operators. In particular, the VEM bilinear forms typically involve two terms. The \emph{polynomial consistency} term acts  on the projection of the discrete functions and is responsible for the accuracy of the method; the \emph{stabilization term} is complementary to the consistency term and is required to ensure  the coercivity of the VEM  bilinear forms. See, for example, \cite{basicprinciples,Eqproj,cangiani2017conforming} for more details including proofs of stability and convergence.

Assume that a discrete solution at  time level $t^n$ has been computed as $\rho_h\in V_h^n$ on the current mesh  $\mathcal{T}_h^n$. When $n=0$, $\rho_h$ is defined as the virtual element interpolation of the initial condition of the PME Problem~\ref{prob::PME}.  Otherwise, $\rho_h$ is the current time discrete solution.

The VEM discretisation of the  velocity potential Problem~\ref{prob::ContinousPotential}  reads: given $\rho_h \in V_h^n$, find $\phi_h \in V_h^n$ such that
\begin{equation} \label{eq::DiscPotential}
    A_h(\phi_h, w_h) = d_h(w_h) \qquad\quad \forall \, w_h\in V_h^n
\end{equation}
with the approximate bilinear form $A_h$ and linear form $d_h$ built by summing element-wise contributions as typical of FEM, hence 
\[A_h(\phi_h, w_h)=\sum_{E\in\mathcal{T}_h^n}A_h^E(\phi_h, w_h)\qquad \text{and}\qquad   d_h(w_h)=\sum_{E\in\mathcal{T}_h^n}d_h^E(w_h).
\] 
As anticipated above however, the preceding definition of the local VEM forms necessitates the use of projections as follows:
\begin{align}
    {A}_h^E(\phi_h,w_h) &= \int_{E}  (\Pi\rho_h)_0
\ \boldsymbol{\Pi}_0 \nabla \phi_h \cdot \boldsymbol{\Pi}_0 \nabla w_h\ d\mathbf{x} + (\Pi\rho_h)_0\, S_A^E(\phi_h-\Pi \phi_h, w_h-\Pi w_h), \label{eq::Atilde_h}\\
    d_h^E(w_h) &= -\int_{E} (\Pi \rho_h)^m_0\ \boldsymbol{\Pi}_0\nabla \rho_h \cdot \boldsymbol{\Pi}_0 \nabla   w_h  \ d\mathbf{x}\label{eq::d_h},
\end{align}
where $ (\Pi\rho_h)_0$ is the constant component of $\Pi\rho_h$. Here,  following~\cite{basicprinciples}, the stabilization form $S_A^E(\cdot,\cdot)$ is defined by
\begin{equation}
    S_A^E(v_h, w_h) = \sum_{l=1}^{m_E} \dof_l(v_h) \cdot \dof_l(w_h),
\end{equation}
with $m_E$ denoting the dimension of $V(E)$, which for the space considered here is equal to the number of vertices of $E$, and with $\text{\dof}_i(w)$ representing the $i$-th degree of freedom of the function $w$. Hence, $\dof_l(v_h)=v_h(\mathbf{x}_l)$ with $\mathbf{x}_l$ denoting the $l$-th vertex of $E$. The integration constant is fixed by constraining a single vertex value of $\phi_h$ to zero. A variety of suitable stabilization choices are admissible \cite{ellipticVEM,cangiani2017conforming} but we adopt the simplest choice in this paper. In particular, in the case when arbitrarily small edges appear in the mesh we refer to \cite{brenner2018virtual} for more appropriate stabilization terms.

The velocity reconstruction Problem~\ref{prob::VelocityReconstruction} is a global $L^2$ projection. Its VEM discretisation reads: given $\phi_h\in V_h^n$, the solution of~\eqref{eq::DiscPotential}, find ${\velocity}_h \in \left[ V_h^n\right]^2$ such that ${\velocity}_h \cdot \mathbf{n} = 0$ on the portion of $\Omega_h^n$ approximating $\partial \Omega_F^n$ and
\begin{equation}
    M_h({\velocity}_h,w_h) = b_h(w_h) \qquad\quad \forall \, w_h \in  V_h^n.
\end{equation}
As before, the forms $M_h$ and $b_h$ are obtained summing the respective elemental forms 
\begin{align}
    M_h^E({\velocity}_h,w_h) &= \int_E \boldsymbol{\Pi} {\velocity}_h\ \Pi w_h\ d\mathbf{x}\ +\ S_M^E({\velocity}_h-\boldsymbol{\Pi}{\velocity}_h, w_h-\boldsymbol{\Pi} w_h), \label{eq::DiscreteMass4Velocity}\\
    b_h^E(w_h) &= \int_E \Pi w_h\ \boldsymbol{\Pi}_0 \nabla \phi_h\ d\mathbf{x},\label{eq::DiscreteRHS4Velocity}
\end{align}
with the stabilization term $S_M^E(\cdot, \cdot)$  given by \cite{Eqproj}
\begin{equation}\label{eq::StabMass}
    S_M^E(\velocity_h,w_h) = |E| \sum_{l=1}^{m_E} \dof_l(\velocity_h) \cdot \dof_l(w_h).
\end{equation}

\subsection{Moving the mesh}
 The mesh is transferred between discrete time levels by displacing the nodes of the mesh and maintaining the mesh connectivity between $t=t^n$ and $t=t^{n+1}$. 
For each mesh node $\mathbf{x}^n$ the new position is obtained by the forward Euler method applied to $\dot{\mathbf{x}}=\velocity_h(\mathbf{x})$, yielding $\mathbf{x}^{n+1}=\mathbf{x}^n+(t^{n+1}-t^n)\velocity_h(\mathbf{x}^n)$.
Thus, consistent with the VEM philosophy,  only the values of $\velocity_h$ at the nodes, that is the degrees of freedom of $\velocity_h$, are required to compute the mesh movement.

\begin{remark}
    We note that the mesh velocities computed in this way do not guarantee a priori that Assumption \ref{assumption::ShapeRegular} holds true indefinitely. For instance, when performing the numerical experiments for contact problems of Section \ref{sec::Tests}, we observed a degradation of mesh quality near the contact boundary. To compensate for this, a harmonic extension operator is used to regularise the velocity field and maintain element sizes. Further details are provided in Section \ref{sec::Tests}. 
\end{remark}

\section{Virtual Element Method for the Solution}\label{sec::VEMSolution}
Once the new mesh node positions have been computed, we consider  the process of updating the solution. 
This is performed in two steps corresponding, respectively, to the conservative ALE update of the mass monitor Problem~\ref{prob::ALEUpdate} and the actual solution update Problem~\ref{prob::SolutionRecon}. 
Before presenting the details on their VEM discretisation, a discussion on the hypothesis leading to such problems is in order. The original moving mesh method in~\cite{MovingReview} was  based on the linear FEM for which the validity at the discrete level  of the material derivatives assumption~\eqref{eq::MaterialDerivative} leading to the ALE update~\eqref{eq::MudotMH} has been proven in~\cite{jimack&Wathen}. In the VEM setting, instead, we exploit the fact that virtual element functions are \emph{only accessed through their nodal values}: in close alignment with the work presented in \cite{dGALEmapping}, only the mesh skeleton velocity is known and used. Hence, in view of the  solution update through the time step $[t^n,t^{n+1}]$, we can assume that ~\eqref{eq::MaterialDerivative} is satisfied by the space-time discrete basis which are then interpolated at the new time level, again, just by accessing the nodal values of the solution.

\subsection{Discretisation of the solution problems}

The initial condition $\rho_h^0$ is approximated by interpolating the degrees of freedom of $\rho^0$ into the VEM space $V_h^0$. Then, the initial mass monitor distribution is computed via
\begin{equation}\label{eq::Mu_h}
    \mu_h^0(w_h^0) = \sum_{E \in \mathcal{T}_h^0} \int_E \Pi \rho_h^0\ \Pi w_h^0\ d\mathbf{x}.
\end{equation}
The next task is the update of the mass monitor over time levels. Once again, this  is performed via the forward Euler method: the new monitor is thus given by $\mu^{n+1}_h(w_h^{n+1}) = \mu^n_h(w_h^{n}) + (t^{n+1} - t^n) \dot{\mu}_h^n(w_h^n)$. This, in turn, requires the approximation of the ALE equation (\ref{eq::ALEUpdateContinuous}) which is performed  still on the old time level (superscript omitted) by
\begin{equation}\label{eq::DiscreteALEUpdate}
    \dot{\mu}_h(w_h) = -\sum_{E \in \mathcal{T}_h^n} \int_E \Pi \rho_h \boldsymbol{\Pi}_0 \nabla w_h \cdot \left\{ \left( \Pi\rho_h \right)_0^{m-1} \boldsymbol{\Pi}_0 \nabla \rho_h + \boldsymbol{\Pi}{\velocity}_h \right\} \ d\mathbf{x}.
\end{equation}

Once the monitor is updated, we set the time level to the new time and update the solution using a VEM discretisation of Problem~\ref{prob::SolutionRecon}: find $\rho_h^{n+1} \in V_h^{n+1}$ such that
\begin{equation}
    m_h(\rho_h^{n+1},w_h^{n+1}) = \mu_h(w_h^{n+1})\ \ \ \ \ \forall \, w_h^{n+1} \in V_h^{n+1}.
\end{equation}

Similarly to the velocity reconstruction, the discrete form $m_h(\cdot,\cdot)$ is computed on the new time level by summing over element contributions
\begin{equation}
    m_h(\rho_h, w_h)=\sum_{E\in\mathcal{T}_h^n}m_h^E(\rho_h, w_h),
\end{equation}
where
\begin{equation}\label{eq::DiscreteMassForm}
    m_h^E(\rho_h,{w}_h) = \int_E \Pi \rho_h\ \Pi w_h\ d\mathbf{x}\ +\ S_m^E(\rho_h, w_h),
\end{equation}
with the stabilization term $S^E_m(\cdot, \cdot)$ being given by
\begin{equation}\label{eq::StabMassScalar}
    S_m^E(\rho_h,w_h) = |E| \sum_{l=1}^{m_E} \dof_l(\rho_h) \cdot \dof_l(w_h).
\end{equation}

\subsection{Partition of unity and conservation}\label{sec::MassConservation}
A beneficial property of the linear virtual element method is that the basis functions form a partition of unity on $\Omega_h^n$ at any discrete time level, i.e. 
\begin{equation}\label{eq::PartitionOfUnity}
    \sum_{i=1}^{\Ndof} \varphi_i = 1,
\end{equation}
where the set $\left\{ \varphi_i \right\}_{i=1}^{\Ndof}$ refers to the set of canonical VEM basis functions associated to the vertices of the mesh \cite{basicprinciples}; that is $\varphi_i(\mathbf{x}_j) = \delta_{ij}$ for $i=1,...,\Ndof$ where $\mathbf{x}_j$ is the j-th node in the mesh.

For a given $\rho_h \in W^n_h$ the monitor integral $\theta^n$ reads
\begin{equation}\label{eq::DiscreteTheta}
    \theta^n = \int_{\Omega^n} \rho_h \ d\mathbf{x},
\end{equation}
from which the polynomial consistency and partition of unity property of the VEM gives
\begin{align}
    \theta^n &= \sum_{E \in \mathcal{T}_h^n} \int_{E} \Pi \rho_h\ d\mathbf{x}\\
    &= \sum_{E \in \mathcal{T}_h^n}  \int_{E} \Pi \rho_h\ \sum_{j=1}^{m_E}\varphi_j\ d\mathbf{x}\\
    &= \sum_{i=1}^{\Ndof} \sum_{E \in \mathcal{T}_h^n} \int_{E} \Pi \rho_h\ \Pi\varphi_i\ d\mathbf{x}.
\end{align}
Therefore the global conservation of the mass monitor is only dependent on the polynomial component of the discrete solution and weighting functions. Further, the exact value of the monitor can also be recovered via
\begin{equation} \label{eq::SumDiscreteMu}
    \theta^n =  \sum_{i=1}^{\Ndof}  \mu_h^n(\varphi_i).
\end{equation}
Finally, considering the partition of unity property, the ALE update equation (\ref{eq::ALEUpdate}) and equation (\ref{eq::SumDiscreteMu}) for the PME, we get
\begin{equation}
    \dot{\theta}^n = 0,
\end{equation}
which agrees with the conservation of mass principle for this particular PDE. In fact, virtual elements preserving relevant global conservation laws in different contexts can be constructed, an example of which is given in \cite{parabolicVEM} for the heat equation.

\section{Implementation Details}\label{sec::Implementation}
This section presents a complete overview of the moving mesh virtual element method. The construction of the required algebraic equations and imposition of boundary conditions are reviewed along with some practical remarks regarding the implementation of this method. The initial weak distribution of the monitor is stored in the vector $\boldsymbol{\mu}^0$ and can be computed using equation (\ref{eq::Mu_h}) whilst the mass matrix $\mathbf{M}^n$ is computed by assembling the contributions from equation (\ref{eq::DiscreteMassForm}).
In order to compute the discrete potential $\phi_h \in V_h^n$ from equation \eqref{eq::DiscPotential}, we solve the linear system $\mathbf{A}^n \boldsymbol{\phi} = \mathbf{d}^n$ with $\mathbf{A}^n$ and $\mathbf{d}^n$ computed using equations \eqref{eq::Atilde_h} and \eqref{eq::d_h}, respectively.
The solution of the resulting linear system determines $\phi_h$ up to an additive constant which is inherited from the continuous formulation of the method; here we impose $\text{\dof}_1(\phi_h) = 0$. Once $\phi_h$ is recovered, the velocity field is reconstructed solving $\mathbf{M}_R^n \mathbf{v} = \mathbf{b}^n$, with $\mathbf{M}_R^n$ and $\mathbf{b}^n$ given by equations \eqref{eq::DiscreteMass4Velocity} and (\ref{eq::DiscreteRHS4Velocity}), respectively.  
The ALE update of vector $\dot{\boldsymbol{\mu}}^{n}$ is then obtained using equation (\ref{eq::DiscreteALEUpdate}) and, along with the mesh nodal velocities $\dot{\mathbf{x}}_i^n = \mathbf{v}_h(\varphi_i)$, provides a system of ODEs which can be approximated using the forward Euler method. The main method is outlined in algorithm \ref{algorithm::MovingMeshMethod}. 

\begin{algorithm} 
\SetAlgoLined
\SetKwInOut{Input}{input}\SetKwInOut{Output}{output}
 \Input{The initial condition $\rho_h^0 \in W_h^0$ and mesh $\mathcal{T}_h^0$, the final time $T$.}

 Set $n=0$\;
 Compute $\boldsymbol{\mu}^0$ according to equation \eqref{eq::Mu_h}\;
\BlankLine
 
 \While{$t^n < T$}{
    Construct and solve $\mathbf{A}^n\boldsymbol{\phi} = \mathbf{d}^n$ using equations (\ref{eq::Atilde_h}) and (\ref{eq::d_h}) for $\phi_h \in W_h^n$\;
    Reconstruct the velocity via $\mathbf{M}_R^n\mathbf{v} = \mathbf{b}^n$\;
    \BlankLine
    Compute the ALE update $\dot{\boldsymbol{\mu}}^n$ from equation (\ref{eq::DiscreteALEUpdate})\;
    \BlankLine
    Select $\Delta t$ and set $t^{n+1} = t^n + \Delta t$\;
    Update the mesh node by $\mathbf{x}^{n+1} = \mathbf{x}^n + \Delta t \mathbf{v}$\;
    Update the monitor distribution by $\boldsymbol{\mu}^{n+1} = \boldsymbol{\mu}^n + \Delta t \dot{\boldsymbol{\mu}}^n$\;
    \BlankLine
    Reconstruct and solve $\mathbf{M}^{n+1} \boldsymbol{\rho}^{n+1} = \boldsymbol{\mu}^{n+1}$ for $\rho_h^{n+1} \in W_h^{n+1}$\;
    Update $n = n + 1$\;
 }
 
 \Output{The final solution $\rho_h^{T}$, the final mesh $\mathcal{T}_{h}^{T}$}
 \caption{Moving mesh VEM}\label{algorithm::MovingMeshMethod}
\end{algorithm}

If required, we can strongly impose any Dirichlet boundary conditions on the solution over time whilst simultaneously conserving the monitor integral, the test space used in recovering the solution is augmented to preserve the partition of unity property (\ref{eq::PartitionOfUnity}). In the solution updates, Dirichlet boundary conditions can be enforced using the methodology presented in \cite{MMFEMBC} which extends to a near identical virtual element approach when considering polygonal elements. This is not a necessary step to attain the accuracy results given in section \ref{sec::Tests} and so the implementation details are not discussed in this paper.

\section{A Contact Algorithm}\label{sec::Contact}
In this section we discuss two occurrences of contact and present corresponding basic node insertion algorithms that allow for localised and minimal changes to the mesh structure. In all mesh refinements considered, the change in mesh topology is only performed at the discrete time-levels.  Hence, to ease notation, the superscript used to denote time-steps is omitted. Modified discrete functions, operators, and vectors are denoted using a hat symbol.

\subsection{Contact scenarios}
Here we present two contact scenarios that are numerically investigated in Section \ref{sec::Tests}. The first scenario concerns the collision of the moving boundary with itself whereas the second situation involves collision with fixed geometric obstacles. 

\begin{figure}
\centering
\begin{minipage}[t]{0.3\textwidth}
\centering
 \includegraphics[scale=0.25]{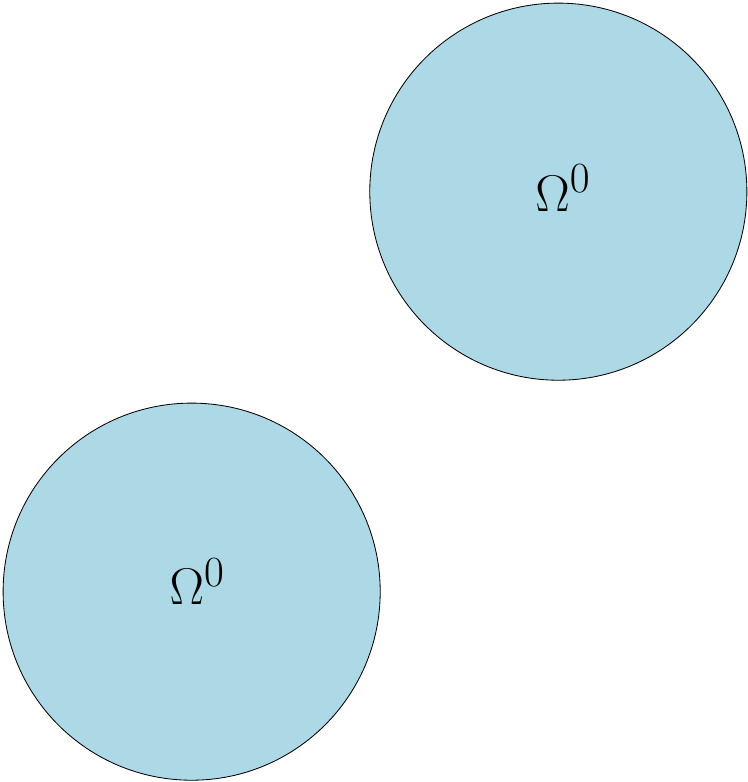}
\end{minipage}
~
\begin{minipage}[t]{0.3\textwidth}
\centering
\includegraphics[scale=0.25]{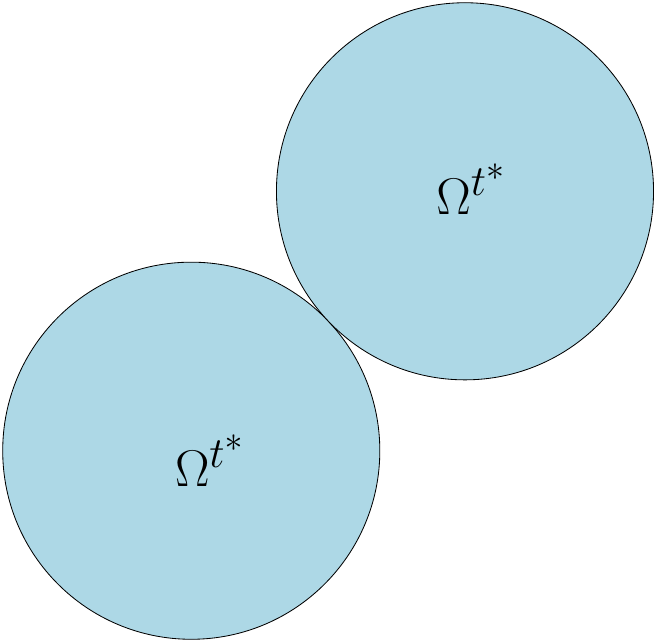}
\end{minipage}
~
\begin{minipage}[t]{0.3\textwidth}
\centering
\includegraphics[scale=0.25]{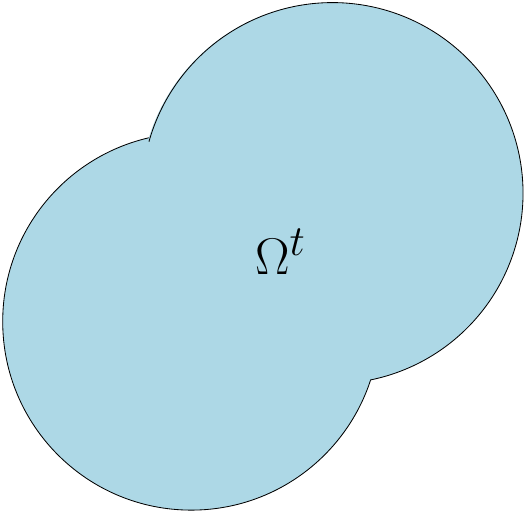}
\end{minipage}
\captionof{figure}{A demonstration of the mesh self-intersection problem. The initial condition has disconnected support in $\tilde{\Omega}$ (left). The disconnected $\partial \Omega^t$ continues to move until some time $t^*$ where the boundary collides with itself (centre). A new connected boundary is formed over time and the topology of $\Omega^t$ is now connected.}
\label{fig::SelfIntersectionExample}
\end{figure}

\emph{Self-intersection} handles the situation where two parts of the moving boundary collide with each other. Typically, a remeshing is required in this instance. By using a VEM, the remeshing can be kept local and simple for colliding elements. A motivational case for a disconnected initial condition of the PME is given in Figure \ref{fig::SelfIntersectionExample}. Moving mesh finite element simulations of this type of problem are presented and discussed in \cite{MovingHessian}.

\begin{figure}
\centering
\begin{minipage}[t]{0.3\textwidth}
\centering
\includegraphics[scale=0.5]{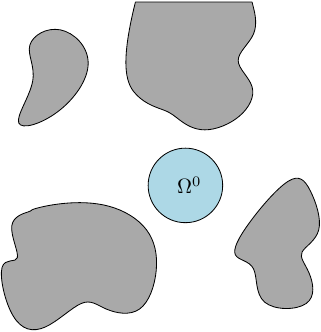}
\end{minipage}
~
\begin{minipage}[t]{0.3\textwidth}
\centering
\includegraphics[scale=0.5]{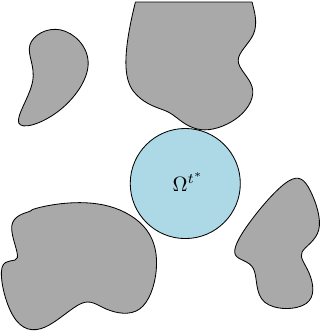}
\end{minipage}
~
\begin{minipage}[t]{0.3\textwidth}
\centering
\includegraphics[scale=0.5]{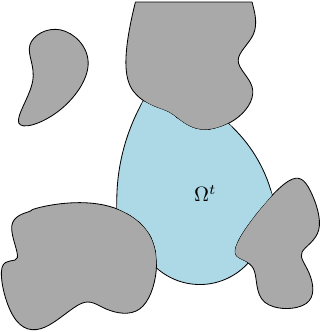}
\end{minipage}
\captionof{figure}{A demonstration of the obstacle contact problem. The initial condition is contained in a region of obstacles (left). After some time $t^*$ the moving boundary collides with the first obstacle (centre). A time-dependent interface now forms between the obstacles and $\partial \Omega^t$ (right).}
\label{fig::ObstacleContactExample}
\end{figure}

\emph{Obstacle contact} is encountered when the evolution of $\Omega^t$ is obstructed by external obstacles. An example of this is the presence of a solid phase in porous media. Figure \ref{fig::ObstacleContactExample} presents an example of collision with impermeable obstacles. By using a collision detection and node insertion algorithm, the moving mesh is capable of simulating the contact and the interaction between the moving mesh and a set of obstacles. As with the self-intersection problem, the VEM allows for this with minimal changes to the mesh topology. Additionally, the VEM is capable of performing local changes to polygonal elements such that the mesh boundary can move around the object boundaries without requiring additional mesh refinements. 

\subsection{Collision detection}
For detecting mesh contact we use an adaptation of the classical \emph{node-to-segment collision detection algorithm} \cite{hallquist1985sliding,wriggers2006computational,VEMcontact}. We only consider boundary mesh and obstacle edges and nodes, thus ensuring that the additional computational cost is $O(N_B^2)$, where $N_B$ is the number of boundary nodes. We consider triplets of points $\left(\mathbf{x}_1^t,\mathbf{x}_2^t,\mathbf{x}_3^t\right)$ where $\mathbf{x}_1^t$ and $\mathbf{x}_2^t$ form a time-dependent edge $e^t$ and $\mathbf{x}_3^t$ is boundary node disconnected from $e^t$. This triplet is referred to as a node-to-edge pair. 

Given a set of linear velocities $\left( \dot{\mathbf{x}}_1, \dot{\mathbf{x}}_2, \dot{\mathbf{x}}_3 \right)$, by defining the vectors
\begin{equation}
    \mathbf{x} = 
    \begin{bmatrix}
        x_1 \\
        x_2 \\
        x_3
    \end{bmatrix}, \ \ 
    \mathbf{y} = 
    \begin{bmatrix}
        y_1 \\
        y_2 \\
        y_3
    \end{bmatrix}, \ \
    \dot{\mathbf{x}} = 
    \begin{bmatrix}
        \dot{x}_1 \\
        \dot{x}_2 \\
        \dot{x}_3
    \end{bmatrix}, \ \
    \dot{\mathbf{y}} = 
    \begin{bmatrix}
        \dot{y}_1 \\
        \dot{y}_2 \\
        \dot{y}_3
    \end{bmatrix}, \ \
\end{equation}
the contact time $\Delta t^*$ between the line of $e^t$ and $\mathbf{x}_3^t$ is given as the minimum positive root of the quadratic equation
\begin{align}
    0 &= a (\Delta t^*)^2 + b (\Delta t^*) + c,\label{eq::ContactQuadratic}
\end{align}
where the coefficients are defined by
\begin{align}
    a = \sum_{i=1}^3 \left( \dot{\mathbf{y}} \times \dot{\mathbf{x}}\right)_i \qquad
    b = \sum_{i=1}^3 \left( \dot{\mathbf{y}} \times {\mathbf{x}} + {\mathbf{y}} \times \dot{\mathbf{x}} \right)_i \qquad
    c = \sum_{i=1}^3 \left( {\mathbf{y}} \times {\mathbf{x}} \right)_i.
\end{align}
In practise, we choose the contact time that makes physical sense (e.g. we discard negative roots as infeasible contact time) and only admit a singular contact time value for $\Delta t^*$. In the case of no feasible contact times we set $\Delta t^* = \infty$ and when two feasible contact times are given by equation \eqref{eq::ContactQuadratic} we choose the minimum of the two. By solving equation (\ref{eq::ContactQuadratic}) for a set of node-to-edge pairing, a set of contact times can be computed. In the context of the moving mesh method, each node-to-edge pair on the boundary of the mesh and (when given) obstacle mesh is considered. If any cases indicate contact, the time step is scaled down to the minimum contact time and the corresponding node-to-edge pair is marked for contact. The detection method is outlined in Algorithm \ref{algorithm::ContactDetection} for a single node-to-edge pairing. 
\begin{algorithm}
\SetAlgoLined
\SetKwInOut{Input}{input}\SetKwInOut{Output}{output}
 \Input{A node-to-edge pairing $\left(\mathbf{x}_{1}^t,\mathbf{x}_{2}^t,\mathbf{x}_{3}^t\right) $, a set of nodal velocities $\left(\dot{\mathbf{x}}_{1}^t,\dot{\mathbf{x}}_{2}^t,\dot{\mathbf{x}}_{3}^t\right)$, the current time step $\Delta t$}

Solve equation \eqref{eq::ContactQuadratic} and set $\Delta t^*$ to the minimum positive root\;
\For{each value of $\Delta t^* \in \mathbb{R}$}{
    Compute $\mathbf{x}_{3}^{t+\Delta t^*}$ and $e^{t+\Delta t^*}$\;
    \uIf{$\Delta t^* \in [0, \Delta t]$ and $\mathbf{x}_{3}^{t+\Delta t^*} \in e^{t+\Delta t^*}$}{
        Mark the node-to-edge pair for contact;
    }
    \Else{
        Set the contact pair to no contact\;
    }
}
 \Output{the node-to-edge contact pair, the contact time step $\Delta t^*$}
 \caption{Contact Detection}\label{algorithm::ContactDetection}
\end{algorithm}

\subsection{Node Insertion algorithm}
Since the mesh allows for general polygonal element shapes, the insertion of a new node into a mesh edge can simply be performed by adding a vertex to the polygons  sharing that edge. Then, a solution value associated with the new vertex must be introduced  which requires an interpolation technique between the old and new global discrete spaces. In \cite{NonHeirVEM}, an elliptic reconstruction operator is employed to preserve the quality of the discrete spatial derivative of the PDE. Instead, here we choose to preserve the polynomial component of the solution $\Pi \rho_h$  between refinements through a redistribution of $\boldsymbol{\mu}^n$.
The reason for this choice is that, by preserving the polynomial component of $\rho_h$, the global mass conservation of $\theta_h^n$ is maintained. This would \emph{not} be the case if interpolation of the degrees of freedom was employed instead. 

On a given element $E \in \mathcal{T}_h$ with an inserted node on the boundary, we impose that the polynomial component of the solution is preserved so that the reconstruction $\hat{\rho}_h \in \hat{V}(E)$ satisfies
\begin{equation}\label{eq::PolynomialLocking}
    \hat{\Pi} \hat{\rho}_h = \Pi \rho_h,
\end{equation}
where $\hat{\Pi}$ denotes the projection operator constructed on the refined VEM space $\hat{V}(E)$. 

The arguments in Section \ref{sec::MassConservation} can be easily modified to show that this approach conserves both locally and globally the mass of the discrete solution. When introducing a new node onto a mesh edge of a given element $E$ under the assumption of equation (\ref{eq::PolynomialLocking}), the local contribution to $\hat{\boldsymbol{\mu}}$ is given by,
\begin{equation}\label{eq::MuUpdate}
    \hat{\mu}_h^{E}(\hat{\varphi}_i) = \int_{E} \Pi \rho_h \hat{\Pi} \hat{\varphi}_i\ d\mathbf{x} \ \ \ i=1,...,\Ndof+1.
\end{equation}
The algorithm for node insertion is given in Algorithm \ref{algorithm::NodeInsertion}.

\begin{algorithm}
\SetAlgoLined
\SetKwInOut{Input}{input}\SetKwInOut{Output}{output}
 \Input{An element $E$, a position $\mathbf{x} \in \partial E$ to insert a node, the monitor distribution $\boldsymbol{\mu}$, the solution $\rho_h \in W_h$, the mesh $\mathcal{T}_{h}$}
 Compute $\Pi \rho_h$ on $E$\;
 Compute the new mesh $\hat{\mathcal{T}}_{h}$ by inserting the node\;
 Compute the new monitor distribution $\hat{\boldsymbol{\mu}}$ using equation (\ref{eq::MuUpdate})\;
 Reconstruct the new solution $\hat{\rho}_h \in \hat{W}_{h}$ by solving $\hat{\mathbf{M}} \hat{\boldsymbol{\rho}} = \hat{\boldsymbol{\mu}}$\;
 \Output{the new solution $\hat{\rho}_h$, the new monitor distribution $\hat{\boldsymbol{\mu}}$, the new mesh $\hat{\mathcal{T}}_{h}$}
 \caption{Node Insertion}\label{algorithm::NodeInsertion}
\end{algorithm}

\subsection{Self-intersection algorithm}
\begin{figure}
\centering
\begin{minipage}[t]{0.3\textwidth}
\centering
\includegraphics[scale=0.5]{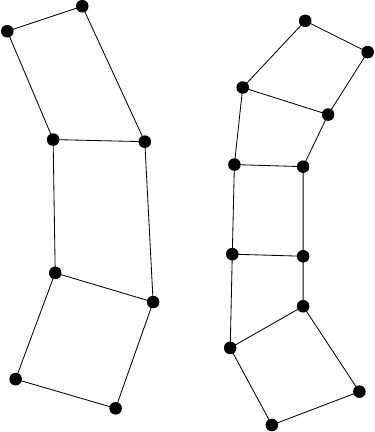}
\end{minipage}
~
\begin{minipage}[t]{0.3\textwidth}
\centering
\includegraphics[scale=0.5]{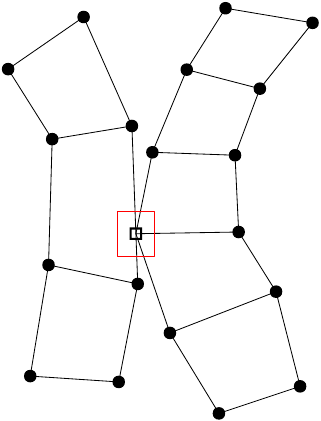}
\end{minipage}
~
\begin{minipage}[t]{0.3\textwidth}
\centering
\includegraphics[scale=0.5]{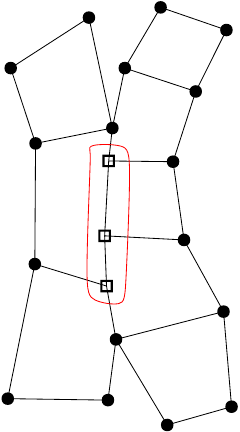}
\end{minipage}
\captionof{figure}{A demonstration of the self intersection algorithm. A sample of elements is shown (left) where a collision is expected. Algorithm \ref{algorithm::Self-intersection} is applied with contact nodes marked by a square (centre). Following subsequent mesh updates certain degrees of freedom no longer lie on the boundary (right). The highlighted degrees of freedom are now treated as internal degrees of freedom.}
\label{fig::SelfIntersectionAlgo}
\end{figure}

Due once more to the VEM  flexibility in element geometries, the self-intersection problem does not require any introduction of additional degrees of freedom. Instead, the local connectivity of the disconnected mesh is updated to include the new node-to-edge pairing. Then, the node insertion algorithm is applied to recompute the solution and monitor distribution. As the boundary node velocities are not arbitrarily set, small edges are likely to appear during node insertion. This has not presented any stability issues within the numerical experiments of section \ref{sec::Tests} and we expect that the method remains robust in the presence of degenerate edges \cite{brenner2018virtual}.
The subsequent mesh velocity problem is then solved based on the updated mesh and corresponding virtual element space. A simple demonstration is provided in Figure \ref{fig::SelfIntersectionAlgo}; the method at a given time step $t^n$ is presented in Algorithm \ref{algorithm::Self-intersection}. 

\begin{algorithm}
\SetAlgoLined
\SetKwInOut{Input}{input}\SetKwInOut{Output}{output}
 \Input{A mesh $\mathcal{T}_h^n$, a solution $\rho_h^n$, a velocity field $\velocity$, a time step size $\Delta t$.}
 Apply Algorithm \ref{algorithm::ContactDetection} for boundary node-to-edge pairs in $\mathcal{T}_h^n$\;
 Update $\Delta t$\;
 Compute $\mathcal{T}_h^{n+1}$, $\boldsymbol{\mu}^{n+1}$, $\rho^{n+1}_h$\;
 \If{any contact pairs are marked}{
    find the element $E \in \mathcal{T}_h^{n+1}$ which contains the marked edge\;
    Apply algorithm \ref{algorithm::NodeInsertion}\;
 }
 Update the boundary conditions of $\rho_h^{n+1}$\;
\caption{Self-intersection}\label{algorithm::Self-intersection}
\end{algorithm}

\subsection{Obstacle contact and pivot node algorithm}
We denote a time-independent polygonal discretization of the set of obstacles by $\mathcal{O}_h$ and consider the mesh node-to-edge pairings  of boundary nodes from $\mathcal{T}_h^n$ and edges of $\mathcal{O}_h$ and vice versa. In the case of obstacle-node to mesh-edge contact, the node insertion Algorithm \ref{algorithm::NodeInsertion} is applied to introduce an additional degree of freedom to the system; we refer to this new node, which is a fixed point on the obstacle geometry, as a ``pivot node''.

For contact between $\mathcal{T}_h^n$ and $\mathcal{O}_h$ a no-penetration condition on the nodal velocities is strongly imposed on the formulation of the potential Problem \ref{prob::ContinousPotential} and velocity reconstruction Problem \ref{prob::VelocityReconstruction}; namely,
\begin{equation}\label{eq::NoPenetration}
    \mathbf{v} \cdot \mathbf{n} = \nabla \phi \cdot \mathbf{n} = 0.
\end{equation}
Hence, movement tangential to the obstacle's boundary is allowed. This is except when a pivot node is introduced, in which case we constrain its velocity to zero to preserve the geometry of the interface between the domain and the obstacle. The obstacle contact algorithm is outlined in Algorithm \ref{algorithm::ObstacleContact}.

Given that the pivot node mesh velocity is constrained to zero, it is possible for the other boundary nodes laying on the obstacle (which have experienced mesh-node to obstacle-edge contact) to pass through the pivot node. When this occurs, the connectivity of the mesh is updated to transfer the pivot node from one mesh edge to another as well as swapping the boundary node from one obstacle face to another. 

Detection for pivot node collision is performed using Algorithm \ref{algorithm::ContactDetection} for connected mesh boundary nodes moving from one obstacle edge to another.

A node is considered to be on $\partial \Omega_{F}^n$ only if both boundary edges sharing that node are in contact with the obstacle. 
We define a node to be ``connected'' if it lies on an edge of the obstacles. Degrees of freedom associated to connected nodes are constrained by equation \eqref{eq::NoPenetration} in the mesh velocity computations. If a node is connected by mesh edges to other connected boundary nodes we consider this to be an ``interface'' node and consequently change the boundary conditions from Dirichlet to Neumann defined in problem \ref{prob::PME} (i.e. we change the degree of freedom from $\partial \Omega_M^n$ to $\partial \Omega_F^n$). If a connected node is not also an interface node, the homogeneous boundary condition and no-penetration condition are maintained. 
The structure of the boundary conditions are updated once every time step.

\begin{algorithm}
\SetAlgoLined
\SetKwInOut{Input}{input}\SetKwInOut{Output}{output}
\Input{A mesh $\mathcal{T}_h^n$, a solution $\rho_h^n$, a velocity field $\velocity$, a time step size $\Delta t$, an obstacle mesh $\mathcal{O}_h$.}

 Apply Algorithm \ref{algorithm::ContactDetection} for boundary nodes in $\mathcal{T}_h^n$ and boundary edges in $\mathcal{O}_h$\;
 Apply Algorithm \ref{algorithm::ContactDetection} for boundary nodes in $\mathcal{O}_h$ and boundary edges in $\mathcal{T}_h^n$\;
 Select the contact pair with the smallest $\Delta t^*$ and set $\Delta t = \Delta t^*$\;
 Compute $\mathcal{T}_h^{n+1}$, $\boldsymbol{\mu}^{n+1}$, $\rho^{n+1}_h$\;
 
 \If{obstacle node to mesh edge is marked}{
    find the element $E \in \mathcal{T}_h^{n+1}$ which contains the marked edge\;
    Apply Algorithm \ref{algorithm::NodeInsertion} to introduce a pivot node $\mathbf{x}_{pivot}$\;
    set $\velocity = \mathbf{0}$ at $\mathbf{x}_{pivot}$\;
 }
 \If{Mesh node to obstacle edge is marked}{
    Set $\velocity \cdot \mathbf{n} = 0$ at mesh node\;
    Update Neumann conditions for $\rho_h^{n+1}$\;
 }
 
\caption{Obstacle contact}\label{algorithm::ObstacleContact}
\end{algorithm}

 When a mesh node $\mathbf{x}_i$ and a pivot node coincide (while the mesh node is moving along the obstacle boundary) the test function associated to a pivot node is chosen to satisfy  $\varphi_{pivot} \equiv \varphi_{i}$. In other words, we duplicate the original VEM basis function and add it to the new discrete space.

\section{Numerical results}\label{sec::Tests}
We report a series of numerical tests for the velocity-based moving mesh virtual element method proposed in the previous sections. Firstly, we present basic convergence test results using a known similarity solution of the PME Problem \ref{prob::PME} for specific choices of velocity and solution recovery. Then, we investigate the effect on the solution of the node insertion algorithms described in Section \ref{sec::Contact}. Finally, we present demonstrations of the contact algorithms of Section~\ref{sec::Contact}.

\subsection{Sample meshes}
\begin{figure}[H]
\centering
\begin{minipage}[t]{.45\textwidth}
  \centering
  \includegraphics[width=1\linewidth]{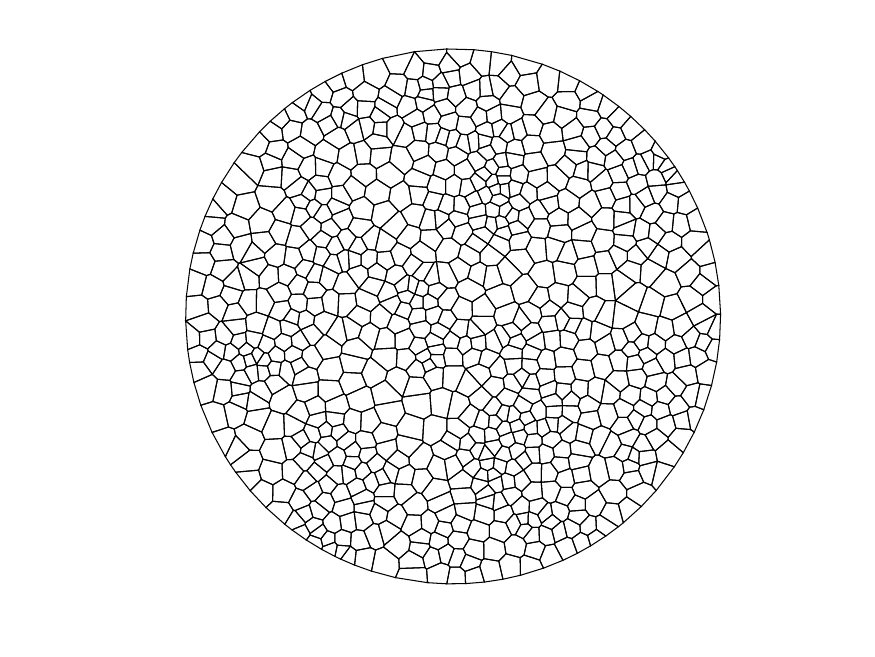}
\end{minipage}
~
\begin{minipage}[t]{.45\textwidth}
  \centering
  \includegraphics[width=1\linewidth]{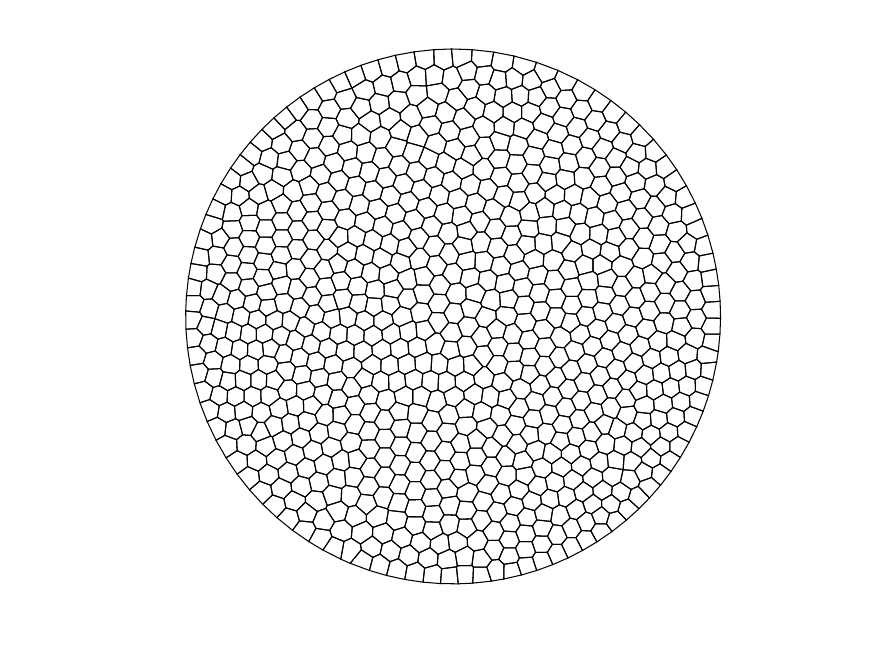}
\end{minipage}
\\
\begin{minipage}[t]{.45\textwidth}
  \centering
  \includegraphics[width=1\linewidth]{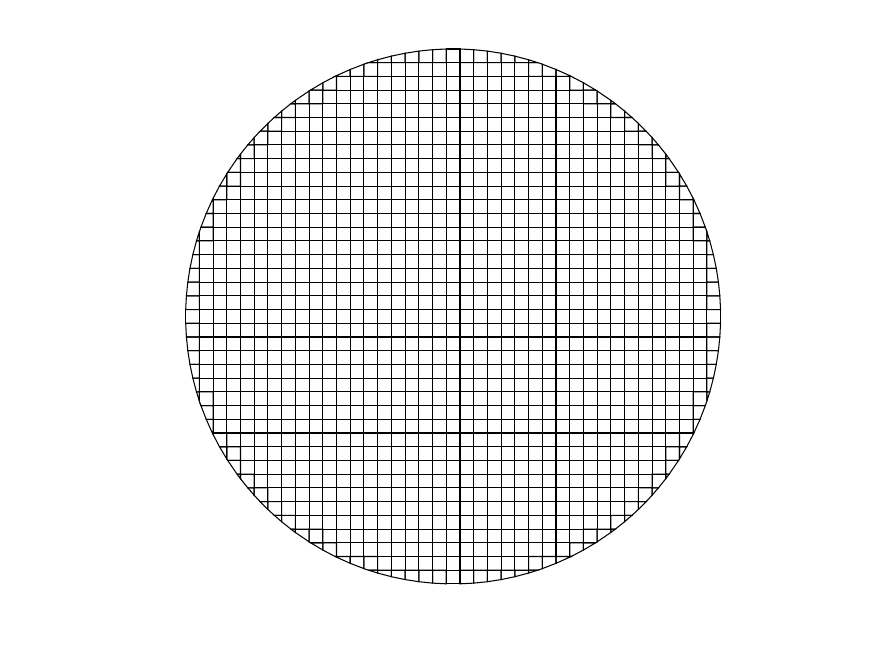}
\end{minipage}
~
\begin{minipage}[t]{.45\textwidth}
  \centering
  \includegraphics[width=1\linewidth]{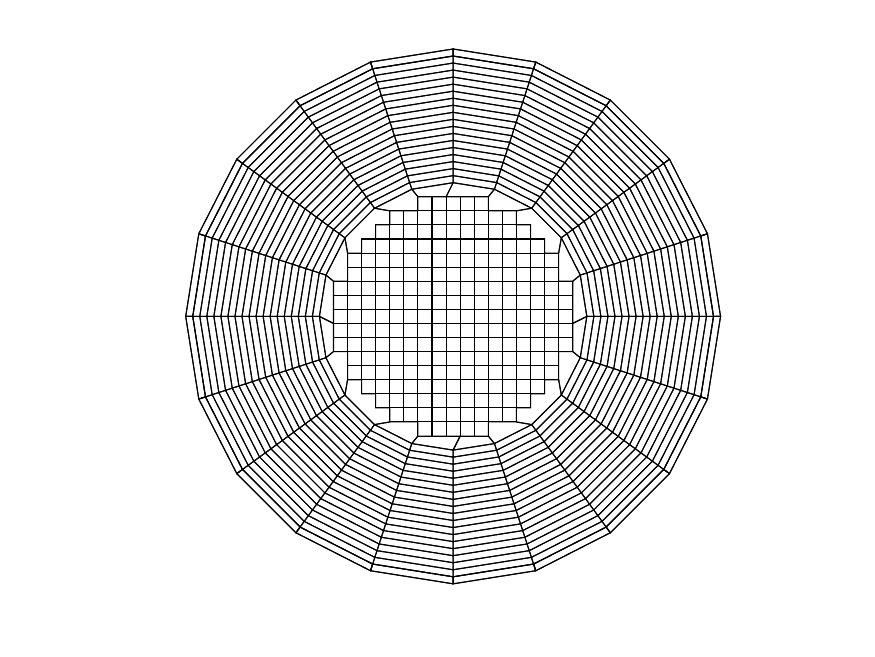}
\end{minipage}

\captionof{figure}{Examples of each of the four mesh types used in numerical tests for a circular domain: the Voronoi Tessellation (top left), the CVT (top right), the grid mesh (bottom left) and the mixed mesh (bottom right).}
\label{fig::MeshTypes}
\end{figure}

We have tested four different mesh types used to subdivide the initial domain; representative examples of each are shown in Figure~\ref{fig::MeshTypes}.
The first mesh is the Voronoi Tessellation produced by randomly sampling mesh seeds in the domain~\cite{Voronoi1,Voronoi2}. The second mesh is a Centroidal Voronoi Tessellation (CVT) produced by the Lloyd algorithm which smooths a given Voronoi tessellation such that the generator points are the barycentric coordinates for each polygon~\cite{lloydconvergence}. The MATLAB package PolyMesher~\cite{polymesher} was used to produce these two mesh types. 
The third mesh is constructed by overlaying the domain with a grid of uniform squares and cutting the mesh along the boundary. The last mesh type is a mixture of uniform Cartesian and polar tessellations. 
Note that the first three initial mesh types may present arbitrarily small edges and, moreover, arbitrarily small elements may appear near the boundary in the grid mesh type, as such potentially contradicting the mesh regularity assumptions stated in Section \ref{sec::MeshAssumptions}.
In this respect, we note that the VEM is known to be quite robust, as we have also witnessed. We refer to the mesh size in each case as the largest element diameter in $\mathcal{T}_h^0$.

\subsection{The PME similarity solution}
There exists a family of radially symmetric solutions on a given initial circular domain of radius $r_0$ for the Problem \ref{prob::PME} defined in \cite{MathBioBook} and given by
\begin{equation}\label{eq::PMEsolution}
\rho(r,t) = 
\begin{cases}
\frac{1}{\lambda(t)^d} \left( 1 - \left( \frac{r}{r_0 \lambda(t)} \right)^2 \right)^{\frac{1}{m}}\ \ \ &|r| \leq r_0\lambda(t) \\
0\ \ \ &\text{otherwise}
\end{cases},
\end{equation}
where $d$ is the spatial dimension, $r_0$ is the initial radius, and
\begin{equation}\nonumber
\lambda(t) = \left( \frac{t}{t_0} \right)^{\frac{1}{2+dm}},\ \ \ t_0 = \frac{r_0^2m}{2(2+dm)}.
\end{equation}
Because of the nature of Equation (\ref{eq::PMEsolution}), the solution
is expected to have finite slope normal to the moving boundary for $m \leq 1$ whilst, for $m > 1$, the solution presents an infinite slope normal to the boundary. Further properties of this analytical solution are discussed in \cite{MovingHessian,MovingReview}.

\subsection{Error computation}
The numerical error is computed for both the solution and mesh by generalising the discrete approximations from \cite{MovingReview}. An $l^1$ solution error is given by,
\begin{equation}\label{eq::SolError}
    \| \rho^n - \rho_h^n \|_{sol} = \frac{1}{\Ndof} \sum_{i=1}^{\Ndof}\ \left| \dof_i(\rho^n) - \dof_i(\rho_h^n) \right|,
\end{equation}
while, for the mesh error, an $l^1$ norm is considered for the radial distance $r(t)$ from the boundary of the mesh to the origin; thus
\begin{equation}\label{eq::MeshError}
    \| r^n - r_h^n \|_{mesh} = \frac{1}{N^B} \sum_{i=1}^{N^B} \left| R_i^n - r_0 \lambda(t^n) \right|.
\end{equation}
Here, $N^B$ denotes the number of boundary nodes in the mesh and $R_i^n$ denotes the radial distance from the origin of the $i$-th boundary node at time $t^n$. A uniform $\Delta t$ that is small enough to ensure numerical stability is set for each initial sample mesh. Note that the meshes used in these numerical tests are not hierarchical. Furthermore, each Voronoi mesh is generated independently from randomly generated seeds. For each reduction of initial mesh size by a factor of 2, the time-step $\Delta t$ is reduced by a factor of 4 to ensure numerical stability. By reducing the time-step size by this factor we also expect that the temporal error to be $O(h^2)$ when using the Forward Euler method for the refinement path of the four mesh types.


\subsection{Convergence test} \label{sec::PMEconvergence}
\begin{figure}
	\centering
	\begin{tikzpicture}
    \begin{groupplot}[group style={
                      group name= myplot,
                      group size= 2 by 4,
                      horizontal sep=1.5cm,
                      vertical sep=1.5cm},
                      height=5.5cm,
                      width=7.15cm,
                      xmode=log,
                      ymode=log,
                      axis background/.style={fill=gray!0}, 
					  legend pos=south east,
					  grid=both,
					  grid style={line width=.1pt, draw=gray!10},
   					  major grid style={line width=.2pt,draw=gray!50}]
                      
        \nextgroupplot
                \addplot+[mark=square, thick, dashed, black, mark options={black, solid}] table [x=Ndof, y=SolError, col sep=comma] {Data Results/voronoiALE.csv};
			    \addplot+[mark=diamond, thick, dashed, blue, mark options={blue, solid}] table [x=Ndof, y=meshError, col sep=comma] {Data Results/voronoiALE.csv};
	    		
	    		\addplot[mark=none, solid, black] coordinates {(1e5,5e-6) (0.25e5,5e-6) (0.25e5,20e-6) (1e5,5e-6)};
	    		\plot[mark=none] (0.25e5,12.5e-6) node[anchor=east] {1};
	    		
	    		\addplot[mark=none, solid, black] coordinates {(1e5,4e-5) (0.25e5,8e-5) (1e5,8e-5) (1e5,4e-5)};
	    		\plot[mark=none] (1e5,12.6e-5) node[anchor=east] {0.5};
		
        \nextgroupplot
                \addplot+[mark=square, thick, dashed, black, mark options={black, solid}] table [x=Ndof, y=SolError, col sep=comma] {Data Results/cvtALE.csv};
			    \addplot+[mark=diamond, thick, dashed, blue, mark options={blue, solid}] table [x=Ndof, y=meshError, col sep=comma] {Data Results/cvtALE.csv};
	    		
	    		\addplot[mark=none, solid, black] coordinates {(1e5,5e-6/2) (0.25e5,5e-6/2) (0.25e5,20e-6/2) (1e5,5e-6/2)};
	    		\plot[mark=none] (0.25e5,12.5e-6/2) node[anchor=east] {1};
	    		
	    		\addplot[mark=none, solid, black] coordinates {(1e5,4e-5) (0.25e5,8e-5) (1e5,8e-5) (1e5,4e-5)};
	    		\plot[mark=none] (1e5,12.6e-5) node[anchor=east] {0.5};
	    		
        \nextgroupplot[xlabel={$\Ndof$}]
                \addplot+[mark=square, thick, dashed, black, mark options={black, solid}] table [x=Ndof, y=SolError, col sep=comma] {Data Results/gridALE.csv};
			    \addplot+[mark=diamond, thick, dashed, blue, mark options={blue, solid}] table [x=Ndof, y=meshError, col sep=comma] {Data Results/gridALE.csv};
	    		
	    		\addplot[mark=none, solid, black] coordinates {(1e5,5e-6/2) (0.25e5,5e-6/2) (0.25e5,20e-6/2) (1e5,5e-6/2)};
	    		\plot[mark=none] (0.25e5,12.5e-6/2) node[anchor=east] {1};
	    		
	    		\addplot[mark=none, solid, black] coordinates {(1e5,6*4e-5/5) (0.25e5,6*8e-5/5) (1e5,6*8e-5/5) (1e5,6*4e-5/5)};
	    		\plot[mark=none] (1e5,6*12.6e-5/5) node[anchor=east] {0.5};
	    		
        \nextgroupplot[xlabel={$\Ndof$}]
                \addplot+[mark=square, thick, dashed, black, mark options={black, solid}] table [x=Ndof, y=SolError, col sep=comma] {Data Results/mixedALE.csv}; \label{plots:SolError}
		    	\addplot+[mark=diamond, thick, dashed, blue, mark options={blue, solid}] table [x=Ndof, y=meshError, col sep=comma] {Data Results/mixedALE.csv}; \label{plots:MeshError}

	    		\addplot[mark=none, solid, black] coordinates {(0.4*1e5,3*5e-6/2) (0.4*0.25e5,3*5e-6/2) (0.4*0.25e5,3*20e-6/2) (0.4*1e5,3*5e-6/2)};
	    		\plot[mark=none] (0.4*0.25e5,3*12.5e-6/2) node[anchor=east] {1};
	    		
	    		\addplot[mark=none, solid, black] coordinates {(0.4*1e5,2*4e-5) (0.4*0.25e5,2*16e-5) (0.4*1e5,2*16e-5) (0.4*1e5,2*4e-5)};
	    		\plot[mark=none] (0.4*1e5,50e-5) node[anchor=east] {1};

    \end{groupplot}
    \path (myplot c1r1.outer north west)
          -- node[anchor=south,rotate=90] {Errors}
          (myplot c1r2.outer south west)
    ;
\path (myplot c1r2.south west|-current bounding box.south)--
      coordinate(legendpos)
      (myplot c2r2.south east|-current bounding box.south);
\matrix[
    matrix of nodes,
    anchor=north,
    draw,
    inner sep=0.4em,
    draw
  ]at([yshift=-1ex]legendpos)
  {
    \ref{plots:SolError}& Solution Error &[8pt]
    \ref{plots:MeshError}& Mesh Error &[8pt] \\
    \\};
\end{tikzpicture}
\caption{PME similarity solution with $m=1$: the $l^1$ solution and mesh errors~\eqref{eq::SolError} and ~\eqref{eq::MeshError}, respectively, at time $T=0.01$ for each mesh type: Voronoi (top left), CVT (top right), Grid (bottom left), Mixed (bottom right).} \label{fig::ALEResults}
\end{figure}
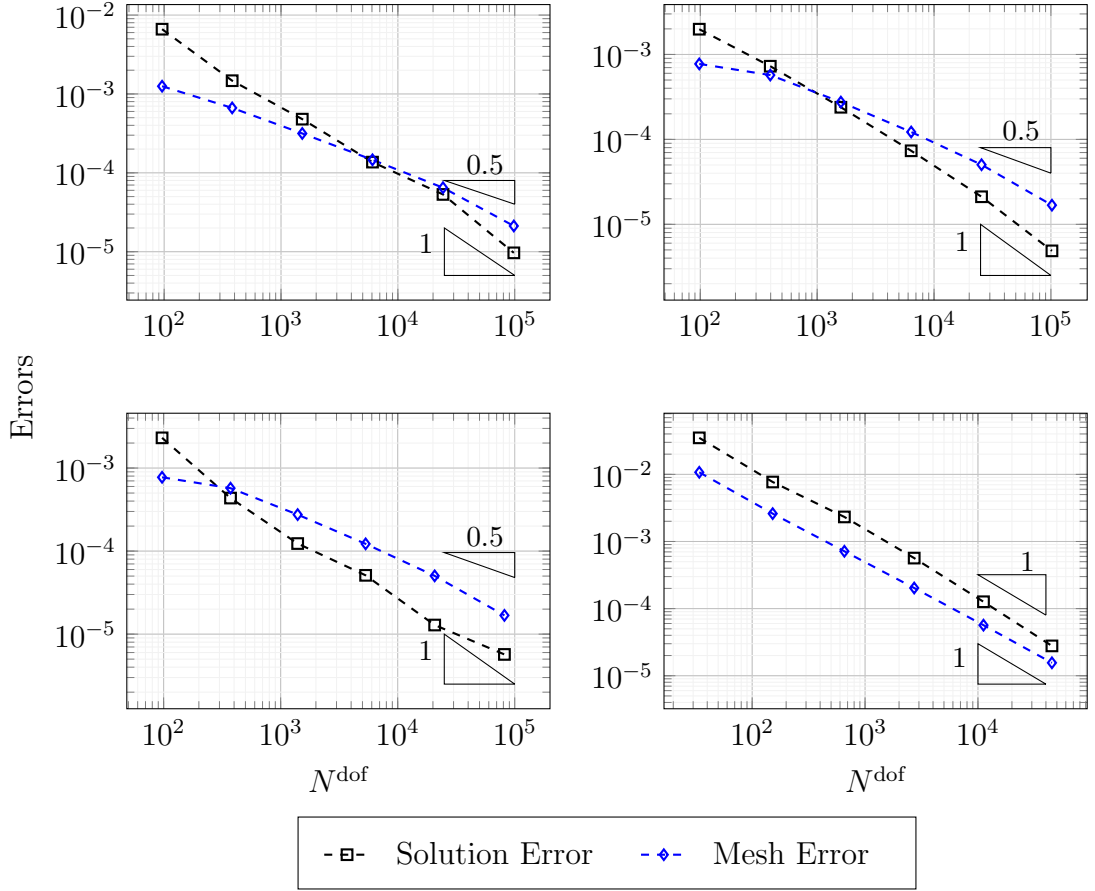

In this first convergence test the solution of the similarity solution for $m=1,d=2$, and $r_0=0.5$ is compared against the numerical solution for $t = t_0 + T$. 
The method is tested on each of the four mesh types for a circular domain. The time step sizes for the coarsest meshes are chosen according to 
\begin{equation}
    \Delta t = \frac{1}{250} h_{mean}^2,    
\end{equation}
where $h_{mean}$ is the average element diameter of the initial mesh $\mathcal{T}_h^0$. In each mesh case we observe that the initial time-step size is approximately $10^{-4}$. From the coarse-mesh time-step sizes we reduce $\Delta t$ by a factor of 4 each time the mesh is refined, which corresponds to the mesh size approximately halving with each refinement.
In choosing the time step sizes we made conservative choices such that the numerical method was stable and presented the expected orders of convergence. A more robust approach would be to use adaptive time-stepping schemes but in this work we present convergence results for a uniform reduction in the time-step.
As shown in Figure~\ref{fig::ALEResults}, second order accuracy is observed for the solution error for all mesh cases when $T=0.01$. 
In the case of the Voronoi and Grid meshes, the empirical order of convergence (EOC) is less smooth compared to the CVT and mixed mesh types. This is most likely due to the weaker shape regularity of elements in these mesh types, but further studies are required.
The mesh error EOC appears to have a long pre-asymptotic regime: the EOC grows monotonically towards the expected rate in all cases with the final computed values ranging between $1.55$ (Voronoi) and $1.83$ (mixed). This is consistent with finite element approximations of Darcy flow which observed the order of convergence of the velocity field to be lower than that of the pressure field \cite{Brezzi_DG_Darcy,Masud_Darcy}.
Conservation of mass in the numerical solution is observed up to machine precision in all test cases.

\begin{remark}
Setting $\dot{\boldsymbol{\mu}} = \mathbf{0}$ produces similar results to those reported in Figure~\ref{fig::ALEResults}. This is referred to as a ``direct recovery'' approach in the literature, see \cite{MovingReview}. For brevity, the corresponding results are not presented  here. 
\end{remark}

\begin{remark}
    When $m>1$, the solution of the PME presents a low regularity as the gradient is unbounded at the moving boundary {\color{red} CITE!}.  Unsurprisingly, applying the moving mesh VEM to the PME with $m=2$, we found that the method remains robust but only attains first order accuracy (results not shown), in line with what already observed for discretizations based on the FEM~\cite{MMFEM}. We expect that appropriately grading the mesh in the vicinity of the moving boundary  can be used to improve the order of convergence as  demonstrated in the FEM case in~\cite{MMFEM}.
\end{remark}

\subsection{Node insertion convergence test}
Our next numerical experiment considers the case of the one-dimensional PME extended in the $x$ direction to a two-dimensional problem. This experiment has two interesting features. Firstly, the initial domain is geometrically exact ($\Omega^0 \equiv \Omega^0_h$) unlike the circular meshes. Secondly, it allows us to test the obstacle contact and node insertion algorithm numerically against a known analytical solution derived from the one-dimensional case of equation \eqref{eq::PMEsolution}. 

This is obtained by considering once again Equation (\ref{eq::PMEsolution}) with the values $m=1,d=1,r_0=0.5$, and $r=y$ on the initial domain given by $\Omega^0 = [-0.5, 0.5]^2$ with initial condition
\begin{equation*}
    \rho(\mathbf{x},0) = 1-4y^2.
\end{equation*}
The mesh is connected to two vertical planes at $x=-0.5$ and $x=0.5$ with a no-penetration condition strongly imposed in the $x$ direction; namely,
\begin{equation*}
    \dot{x} = 0 \ \ \ \text{when } |x| = 1/2.
\end{equation*}

\begin{figure}
\centering
\begin{minipage}[t]{.45\textwidth}
  \centering
  \includegraphics[width=1\linewidth]{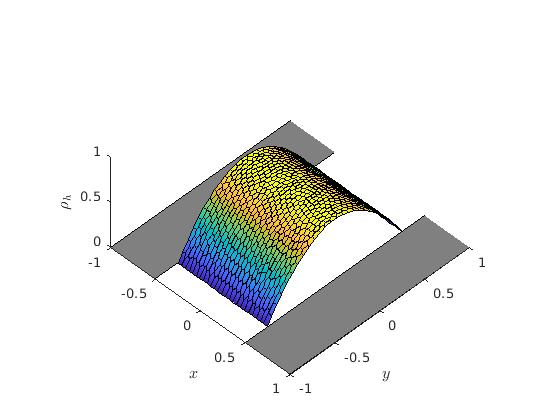}
\end{minipage}
~
\begin{minipage}[t]{.45\textwidth}
  \centering
  \includegraphics[width=1\linewidth]{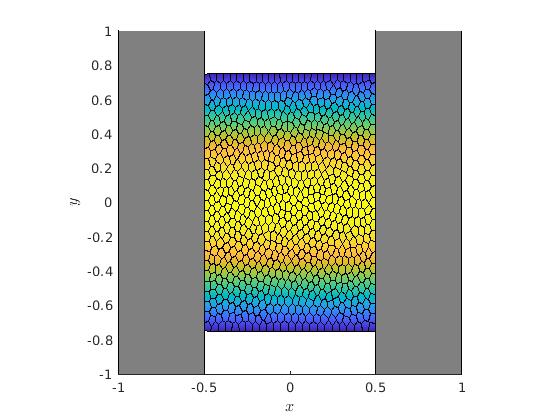}
\end{minipage}
\captionof{figure}{A solution snapshot at time $T=0.1$ for a CVT mesh with $800$ elements (right). }
\label{fig::1DResult}
\end{figure}

Solution snapshots at time $T = 0.1$  are shown in Figure \ref{fig::1DResult} for the CVT mesh type. The mesh error is exclusively computed on the top and bottom faces of the rectangular domain. 

\begin{figure}
    \centering
\begin{subfigure}[b]{0.55\textwidth}
	\centering
	\begin{tikzpicture}
		\begin{axis}[xmode=log,
					 ymode=log,
					 xlabel=$\Ndof$,
					 ylabel=$Errors$, 
					 width=\textwidth, 
					 axis background/.style={fill=gray!0}, 
					 legend pos=south west,
					 grid=both,
					 grid style={line width=.1pt, draw=gray!10},
   					 major grid style={line width=.2pt,draw=gray!50}]
			\addplot+[mark=square, thick, dashed, black, mark options={black, solid}] table [x=Ndof, y=SolError, col sep=comma] {Data Results/pivotCVT.csv};
			\addplot+[mark=diamond, thick, dashed, blue, mark options={blue, solid}] table [x=Ndof, y=meshError, col sep=comma] {Data Results/pivotCVT.csv};
					
			\addplot[mark=none, solid, black] coordinates {(0.2*1e5,5*3*5e-6/2) (0.2*0.25e5,5*3*5e-6/2) (0.2*0.25e5,5*3*20e-6/2) (0.2*1e5,5*3*5e-6/2)};
	    		\plot[mark=none] (0.2*0.25e5,5*3*12.5e-6/2) node[anchor=east] {1};
	    		
	    	\addplot[mark=none, solid, black] coordinates {(0.2*1e5,5*2*4e-5) (0.2*0.25e5,5*2*16e-5) (0.2*1e5,5*2*16e-5) (0.2*1e5,5*2*4e-5)};
	    	\plot[mark=none] (0.2*1e5,5*50e-5) node[anchor=east] {1};
		\end{axis}
		
\path (myplot c1r1.south west|-current bounding box.south)--
      coordinate(legendpos)
      (myplot c1r1.south east|-current bounding box.south);
\matrix[
    matrix of nodes,
    anchor=north,
    draw,
    inner sep=0.4em,
    draw
  ]at([yshift=-1ex]legendpos)
  {
    \ref{plots:SolError}& Solution Error &[8pt]
    \ref{plots:MeshError}& Mesh Error &[8pt] \\
    \\};
	\end{tikzpicture}
\end{subfigure}
\label{fig::PivotTests}
\caption{Node insertion convergence test on a 1D-type PME similarity solution with $m=1$: the $l^1$ solution and mesh errors~\eqref{eq::SolError} and ~\eqref{eq::MeshError}, respectively, at time $T=0.1$.
}\label{fig::PivotResults}
\end{figure}
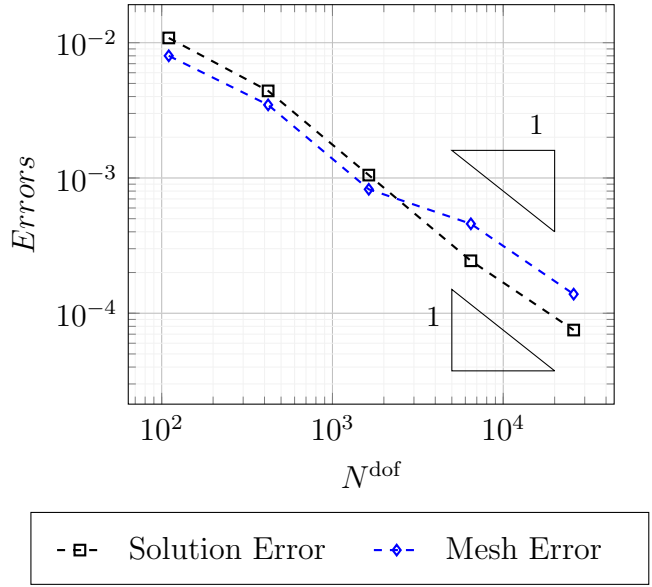

In this test we only focus our attention on the CVT mesh type. 
We test the accuracy of the node insertion algorithm by including a discretization of the two planes into intervals in the $y$ direction. In reference to the fixed domain PME \ref{eq:PME1} we define $\tilde{\Omega} := [-1,1]^2$ and discretize the boundary $\partial \tilde{\Omega}$ into N uniformly spaced intervals to construct an $N$-gon. Vertices are then removed that intersect $\Omega^0$. This results in uniformly discretized intervals along both contact planes.

 For the discretization of the contact planes $N$ is set to an initial value of $32$ and is doubled with each mesh refinement. Convergence results are reported in Figure~\ref{fig::PivotResults}. Here we observe  second order accuracy in both the solution and mesh error. Also the mass is conserved, by design, up to machine precision throughout each test. 

\subsection{Contact demonstrations}
We finally present two demonstrations of the node insertion algorithms for Problem~\ref{prob::PME} in challenging scenarios without known analytical solution.

To ensure the quality of the mesh is preserved we introduce an alternative velocity field on the interior of the moving mesh based on the ALE approach. First, the method outlined in Section \ref{sec::Implementation} is applied to approximate the Lagrangian boundary velocity $\mathbf{v}$. The mesh velocity $\mathbf{v}$ for the internal node movement is then replaced by $\tilde{\velocity}$, the harmonic extension of the Lagrangian boundary velocity: as such it is the solution to the problem
\begin{align}
    -\Delta \tilde{\mathbf{v}} &= \mathbf{0}\\
    \tilde{\velocity}|_{\partial \Omega^t} &= \mathbf{v}.
\end{align}
In a VEM framework this alternative mesh velocity is computed by solving a standard Poisson's equation with Dirichlet boundary conditions \cite{hitchiker}. Other choices for the interior mesh velocities include a modified monitor function $\mathbb{M}(\rho)$ \cite{MMFEMmonitors}, pseudo-elastic, and biharmonic formulations \cite{Richter2017}. We remark that the quality of the mesh will still deteriorate over time. The purpose of these examples is to demonstrate the application of the node insertion algorithms. Optimising the choice of ALE velocity is left for future investigation.

\begin{figure}
\centering
\begin{minipage}[t]{.45\textwidth}
  \centering
  \includegraphics[width=1\linewidth]{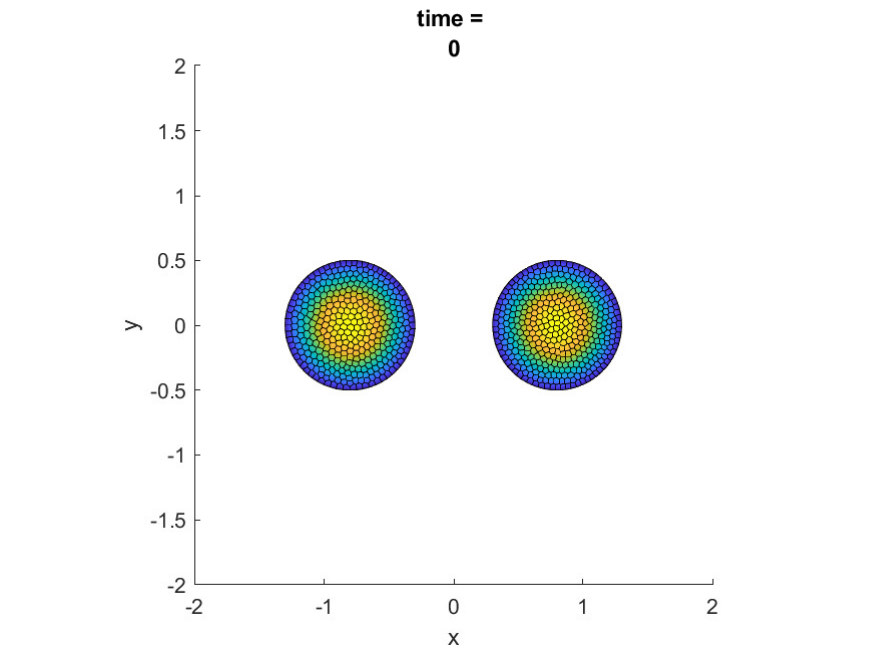}
\end{minipage}
~
\begin{minipage}[t]{.45\textwidth}
  \centering
  \includegraphics[width=1\linewidth]{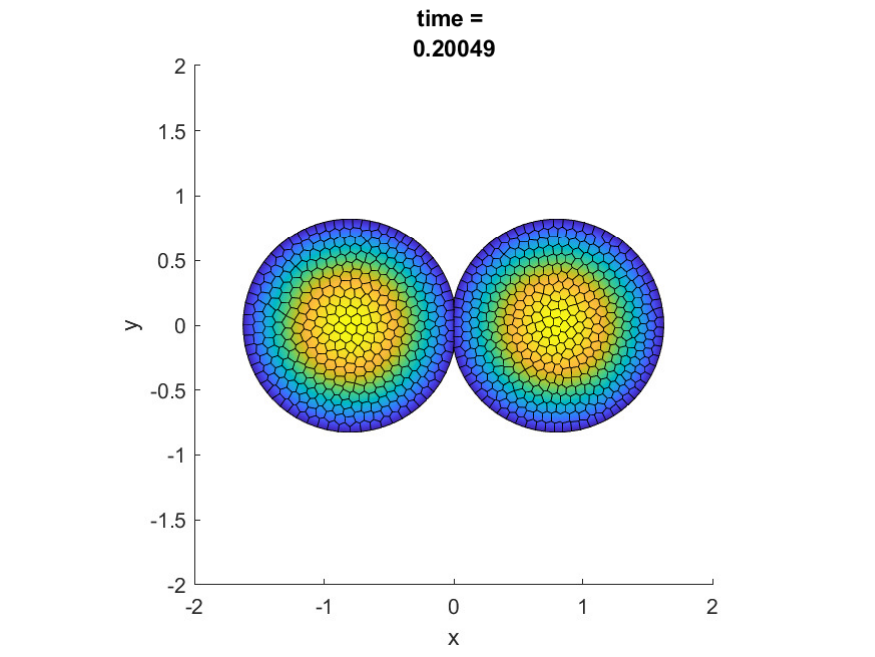}
\end{minipage}
\\
\begin{minipage}[t]{.45\textwidth}
  \centering
  \includegraphics[width=1\linewidth]{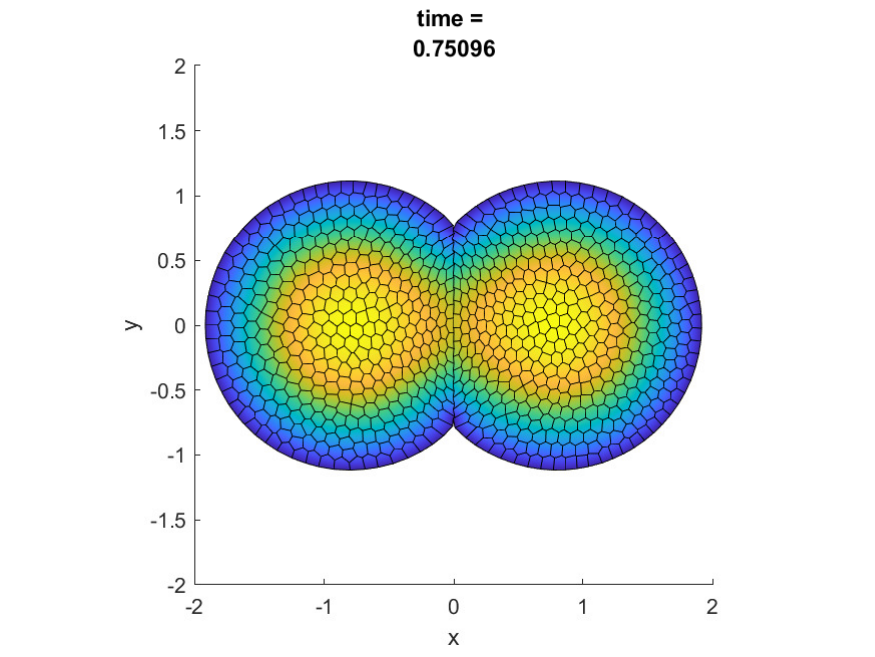}
\end{minipage}
~
\begin{minipage}[t]{.45\textwidth}
  \centering
  \includegraphics[width=1\linewidth]{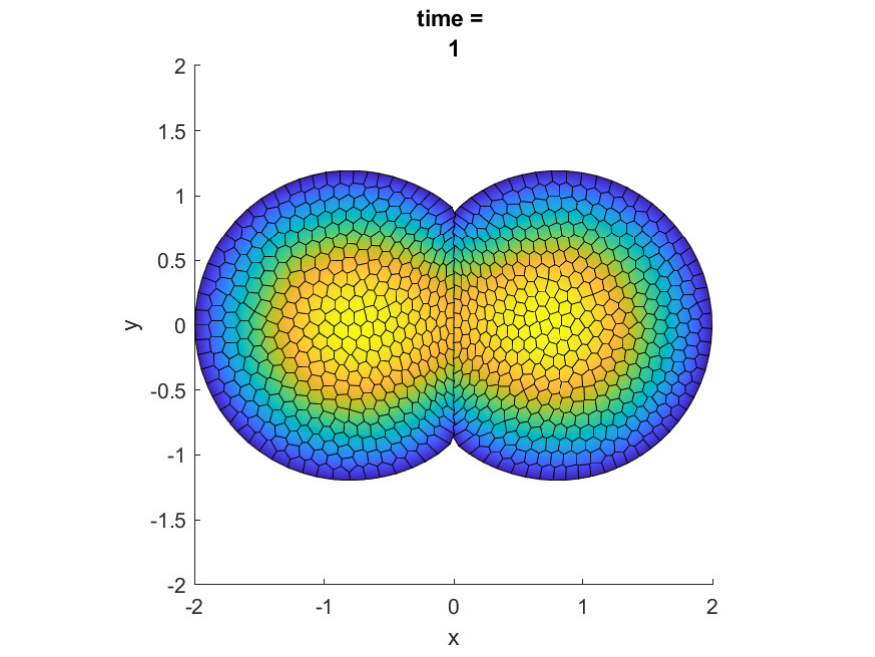}
\end{minipage}
\captionof{figure}{Self-intersection demonstration: snapshots of the moving mesh VEM  solution at $t=0$ (top left), $t=0.20049$ (top right), $t=0.75096$ (bottom left), and $t=1$ (bottom right). A  CVT type mesh  with 800 elements was used to initialise the mesh at $t=0$.}
\label{fig::TwoMesh}
\end{figure}

In the first demonstration, we consider an initial condition of the PME that has a disconnected support such that self-intersection is expected to occur. The initial condition is given by
\begin{equation}
    \rho(\mathbf{x},0) = 
\begin{cases}
1-4r_1^2\ \ \ &r_1 = | \mathbf{x}-(-0.8,0)|,\ r_1 \leq 1/2, \\
1-4r_2^2\ \ \ &r_2 = | \mathbf{x}-(0.8,0)|,\ r_2 \leq 1/2, \\
0\ \ \ &\text{otherwise}.
\end{cases}
\end{equation}
An illustrative example of such initial condition is given in Figure~\ref{fig::TwoMesh} (top-left plot). The standard method is applied to simulate the PME for $m=1$ with the contact detection Algorithm~\ref{algorithm::ContactDetection} applied at every time level to check for collision between elements. When contact occurs, Algorithm~\ref{algorithm::Self-intersection} is used to update the monitor distribution whilst the Dirichlet boundary degrees of freedom are flagged as interior degrees of freedom as the mesh connectivity evolves. Snapshots of the solution evolving over time are reported in  Figure~\ref{fig::TwoMesh}. The behaviour of the PME solution over time is in agreement with fixed mesh finite element approximations of this problem and similar benchmark tests performed for the PME in \cite{MovingHessian}.

\begin{figure}
\centering
\begin{minipage}[t]{.45\textwidth}
  \centering
  \includegraphics[width=1\linewidth]{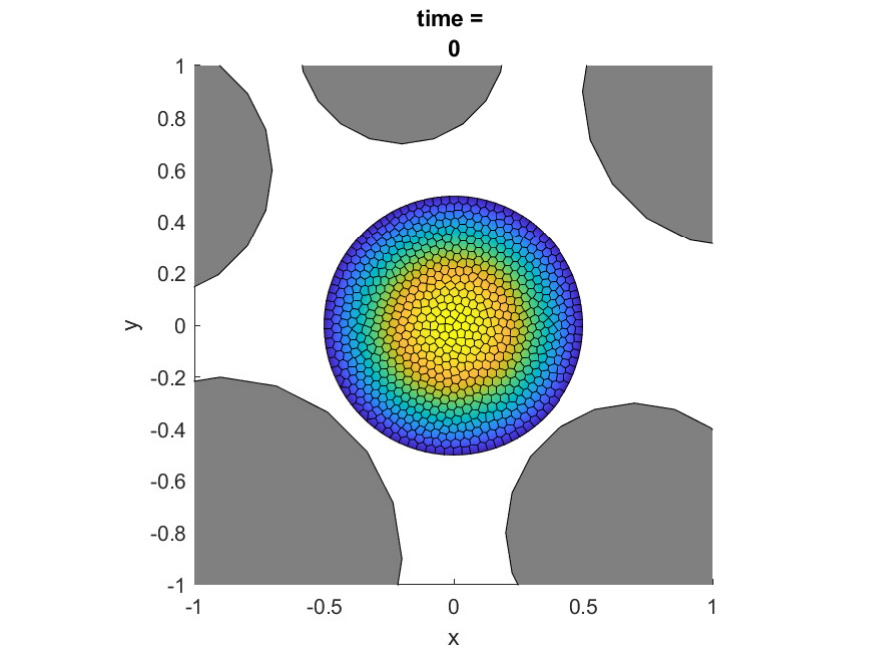}
\end{minipage}
~
\begin{minipage}[t]{.45\textwidth}
  \centering
  \includegraphics[width=1\linewidth]{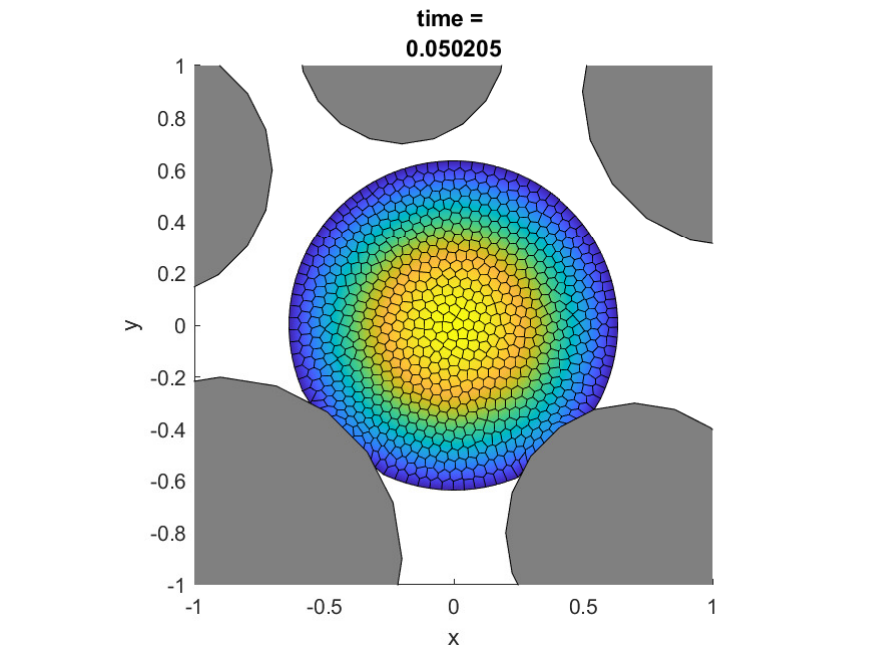}
\end{minipage}
\\
\begin{minipage}[t]{.45\textwidth}
  \centering
  \includegraphics[width=1\linewidth]{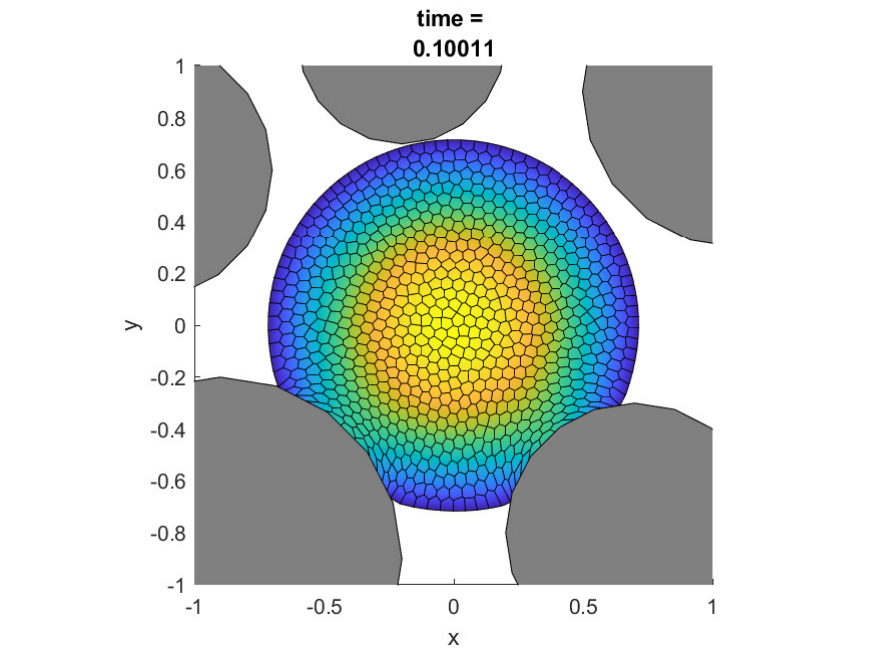}
\end{minipage}
~
\begin{minipage}[t]{.45\textwidth}
  \centering
  \includegraphics[width=1\linewidth]{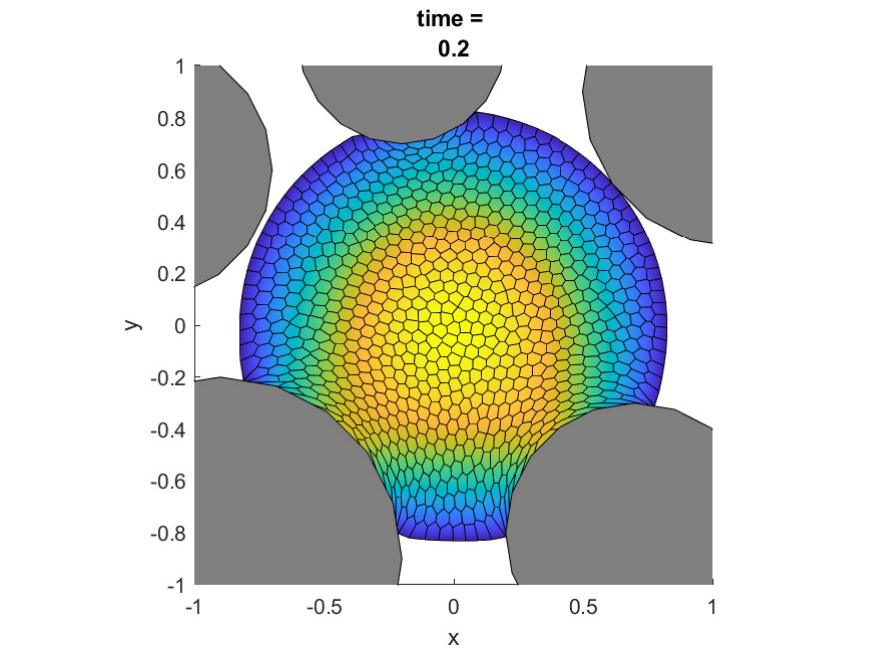}
\end{minipage}
\captionof{figure}{Obstacle contact demonstration: snapshots of the moving mesh VEM solution at $T=0$ (top left), $T=0.050205$ (top right), $T=0.10011$ (bottom left) and $T=0.2$ (bottom right). A  CVT type mesh  with 800 elements was used to initialise the mesh at $t=0$.}
\label{fig::SwissCheese}
\end{figure}

To demonstrate the obstacle contact algorithm we consider again the initial condition  given by Equation (\ref{eq::PMEsolution}). A set of obstacles are added to the computational domain in the shape of circles with random radii and centres. Each circle is approximated as a uniform polygon of comparable accuracy to the initial mesh. We tested the moving mesh VEM starting with a CVT type mesh made of 800 elements discretising the support of the initial solution.
A  few snapshots of the numerical solution are shown in Figure~\ref{fig::SwissCheese}.  Pivot nodes are inserted and removed along the contact. Mesh degeneration occurs after $T=0.2$ and thus the simulation had to be terminated. In this case, and similar to finite element methods, a remeshing approach would rectify this issue. 

\subsection{A Fourth-Order Problem}
To demonstrate the extensibility of the moving mesh VEM we consider the following fourth-order nonlinear diffusion problem used as a benchmark for the original moving mesh method \cite{MovingReview,MMFEM}, for which the whole of the boundary $\partial \Omega^t$ is free to move, and the differential operator is given by $\mathcal{L} = -\nabla \cdot (\rho^m \nabla \Delta \rho)$. In this work we choose $m = 1$, for which there is a simple similarity solution, defined below. The resulting time-dependent equation $\partial\rho/\partial t =\mathcal{L}$ is complemented, at the free boundary, with two conditions on $\rho$, namely $ \rho = \nabla \rho \cdot \mathbf{n} =  0$ plus the kinematic condition
$\rho \mathbf{v} \cdot \mathbf{n} = \rho \nabla \Delta \rho \cdot \mathbf{n}$, which is used to determine the boundary velocity $\mathbf{v}$.

In view of its numerical solution, we re-write the fourth-order problem as a coupled system of second-order PDEs by introducing a pressure term $p = -\Delta \rho$. The problem then reads as: \emph{find $\rho = \rho(\mathbf{x},t)$ such that $\rho(\mathbf{x},0) = \rho^0(\mathbf{x})$ for $\mathbf{x} \in \Omega^0$ and, for all $t \in (0,T]$,}
\begin{align}
    \frac{\partial \rho}{\partial t} &= \nabla \cdot (\rho \nabla p) \qquad &&\mathbf{x} \in \Omega^t,\\
    \label{eq:peq}
    p &= - \Delta \rho \qquad &&\mathbf{x} \in \Omega^t, \\
    \rho = \nabla \rho \cdot \mathbf{n} &=  0 \qquad &&\mathbf{x} \in \partial \Omega^t, \\
    \rho \mathbf{v} \cdot \mathbf{n} &= - \rho \nabla p \cdot \mathbf{n} \qquad &&\mathbf{x} \in \partial \Omega^t.
\end{align}
This problem is structurally very similar to the porous medium equation problem and the moving mesh algorithm remains mostly unchanged. The main addition is an intermediate step which provides the pressure by discretising the weak form of~\eqref{eq:peq}, namely: \emph{given $\rho \in H^1(\Omega^t)$, find $p \in H^1(\Omega^t)$ such that}
\begin{align}
    \int_{\Omega^t} pw\ d\mathbf{x} = \int_{\Omega^t} \nabla \rho \cdot \nabla w\ d\mathbf{x} \qquad \forall w \in H^1(\Omega^t).
\end{align}
This is discretised using once again the VEM applied to the problem within each time-step $t^n$: given $\rho_h \in V_h^n$ find $p_h \in V_h^n$ such that
\begin{align}
    m_h(p_h, w_h) = \sum_{E \in \mathcal{T}^n_h} \int_E \boldsymbol{\Pi}_0 \nabla \rho_h \cdot \boldsymbol{\Pi}_0 \nabla w_h\ d\mathbf{x} \qquad \forall w_h \in V_h^n,
\end{align}
where $m_h(\cdot,\cdot)$ is defined by Equations \eqref{eq::DiscreteMassForm} and \eqref{eq::StabMassScalar}.

In Problems \ref{prob::ContinousPotential} and \ref{prob::ALEUpdate} the right-hand side terms are also modified for this problem to give, respectively, 
\begin{align}
    d(w) &= -\int_{\Omega^t} \rho \nabla p \cdot \nabla w\ d\mathbf{x}, \\
    \dot{\mu}^t(w) &= \int_{\omegat} - \rho \nabla w \cdot \left\{ \nabla p + {\velocity} \right\} \ d\mathbf{x} \qquad\quad \forall \, w \in H^1(\omegat).
\end{align}
These equations are approximated using the VEM discretizations \eqref{eq::d_h} and \eqref{eq::DiscreteALEUpdate}, as described in Sections \ref{sec::VEMVelocity} and \ref{sec::VEMSolution} for the approximation of the corresponding integrals for the porous medium equation. These are computed at time $t^n$ as follows:
\begin{align}
    d_h(w_h) &= -\sum_{E \in \mathcal{T}_h^n} \int_{E} (\Pi \rho_h)_0\ \boldsymbol{\Pi}_0\nabla p_h \cdot \boldsymbol{\Pi}_0 \nabla   w_h  \ d\mathbf{x}, \\
    \dot{\mu}_h(w_h) &= -\sum_{E \in \mathcal{T}_h^n} \int_E \Pi \rho_h \boldsymbol{\Pi}_0 \nabla w_h \cdot \left\{\boldsymbol{\Pi}_0 \nabla p_h + \boldsymbol{\Pi}{\velocity}_h \right\} \ d\mathbf{x},
\end{align}
$\forall w_h \in V_h^n$.

In view of assessing numerically the resulting moving mesh VEM, we recall  that this fourth-order nonlinear diffusion problem has a radially symmetric similarity solution given by
\begin{equation}\label{eq::FourthOrdersolution}
\rho(r,t) = 
\begin{cases}
At^{\beta}U_0 (1-\eta^2)^2 &|r| \leq A^{\frac{1}{4}} t^{\delta} \\
0\ \ \ &\text{otherwise}
\end{cases},
\end{equation}
where 
\begin{equation}
    \eta = \frac{r}{A^{\frac{1}{4}}t^\delta}, \qquad \delta = \frac{1}{4+d}, \qquad \beta = 4\delta -1, \qquad A=U_0^{-4\delta}.
\end{equation}
Setting $d=2$, we fix  $U_0 = 1/192$ so that $\rho(0,t_0) = 1$ and $t_0 = 1/192$ is specified so that the initial radius is equal to $1$.

\begin{figure}
    \centering
\begin{subfigure}[b]{0.55\textwidth}
	\centering
	\begin{tikzpicture}
		\begin{axis}[xmode=log,
					 ymode=log,
					 xlabel=$\Ndof$,
					 ylabel=$Errors$, 
					 width=\textwidth, 
					 axis background/.style={fill=gray!0}, 
					 legend pos=south west,
					 grid=both,
					 grid style={line width=.1pt, draw=gray!10},
   					 major grid style={line width=.2pt,draw=gray!50}]
			\addplot+[mark=square, thick, dashed, black, mark options={black, solid}] table [x=Ndof, y=SolError, col sep=comma] {Data Results/fourthOrder.csv};
			\addplot+[mark=diamond, thick, dashed, blue, mark options={blue, solid}] table [x=Ndof, y=meshError, col sep=comma] {Data Results/fourthOrder.csv};
					
			\addplot[mark=none, solid, black] coordinates {(0.1*1e5,4*5*3*5e-6/2) (0.1*0.25e5,4*5*3*5e-6/2) (0.1*0.25e5,4*5*3*20e-6/2) (0.1*1e5,4*5*3*5e-6/2)};
	    		\plot[mark=none] (0.1*0.25e5,4*5*3*12.5e-6/2) node[anchor=east] {1};
	    		
	    	\addplot[mark=none, solid, black] coordinates {(0.1*1e5,2*5*2*4e-5) (0.1*0.25e5,2*5*2*16e-5) (0.1*1e5,2*5*2*16e-5) (0.1*1e5,2*5*2*4e-5)};
	    	\plot[mark=none] (0.1*1e5,2*5*50e-5) node[anchor=east] {1};
		\end{axis}
		
\path (myplot c1r1.south west|-current bounding box.south)--
      coordinate(legendpos)
      (myplot c1r1.south east|-current bounding box.south);
\matrix[
    matrix of nodes,
    anchor=north,
    draw,
    inner sep=0.4em,
    draw
  ]at([yshift=-1ex]legendpos)
  {
    \ref{plots:SolError}& Solution Error &[8pt]
    \ref{plots:MeshError}& Mesh Error &[8pt] \\
    \\};
	\end{tikzpicture}
\end{subfigure}
\caption{Convergence test for the fourth-order diffusion problem: the $l^1$ solution and mesh errors~\eqref{eq::SolError} and ~\eqref{eq::MeshError}, respectively, at time $T=0.01$ using the CVT mesh type.
}\label{fig::FourthOrderResults}
\end{figure}
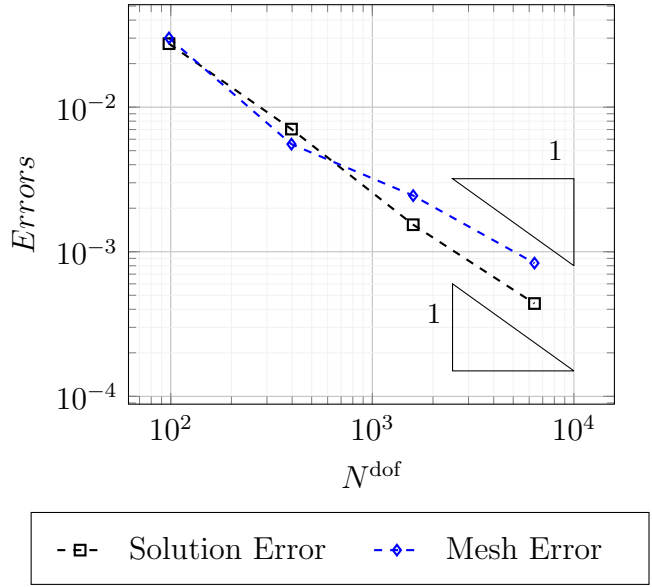

The VEM is tested on the same sequence of CVT-type meshes used in Section \ref{sec::PMEconvergence}, with the same coarse-mesh time-step size of $10^{-4}$ and a reduction by a factor of $4$ each time the mesh is refined. Figure \ref{fig::FourthOrderResults} shows that second-order accuracy is again attained for both the solution and mesh errors. Similar to the PME, the mass is conserved exactly at each time step.

\section{Conclusions}\label{sec::conclusions}
In this paper we have combined, for the first time, a velocity-based moving mesh method with a virtual element method. This was achieved by extending the FEM formulation of the moving mesh method to a linear virtual element scheme on polygonal meshes. Numerical tests for the porous medium equation and a fourth-order diffusion problem show that the proposed method obtains the same orders of accuracy as the original finite element approach. In fact, given that the linear VEM reduces to the linear FEM on triangular elements, our approach provides a natural extension of moving mesh FEM to polygonal mesh settings. Demonstrations of node insertion algorithms suggest that the virtual element method offers practical extensions to more complex problems. In particular, this work shows that it is  straightforward to deal with situations where the moving boundary meets fixed obstacles or merges with other moving boundaries.

The VEM approach provides a flexible discretisation framework for adaptive moving mesh methods. For instance, VEM with curved elements are being developed with a level of generality out of reach for standard FEM, see eg.~\cite{Polycurved,dGALEmapping,Artioli,Ferguson22}.
As such, the VEM is more suitable for the generalisation of moving mesh approaches in the higher-order setting, including for the solution of challenging problems such as phase-field, fluid-structure interaction, and moving surface PDEs~\cite{MovingReview},
including in three-dimensions~\cite{MFDcurved06,dGALEmapping,Bend3d}. We remark that these are largely open problems also for the more standard moving mesh finite element method.
Extensions of the moving mesh VEM in these directions will be the subject of future works.

\section*{Acknowledgements}
The corresponding author was supported by an EPSRC doctoral training grant.


\begin{thebibliography}{10}

\bibitem{Eqproj}
B.~Ahmad, A.~Alsaedi, F.~Brezzi, L.~D. Marini, and A.~Russo.
\newblock Equivalent projectors for virtual element methods.
\newblock {\em Computers \& Mathematics with Applications}, 66(3):376--391,
  2013.

\bibitem{FSIPolyDG}
P.F. Antonietti, M.~Verani, C.~Vergara, and S.~Zonca.
\newblock Numerical solution of fluid-structure interaction problems by means
  of a high order discontinuous galerkin method on polygonal grids.
\newblock {\em Finite Elements in Analysis and Design}, 159:1--14, 2019.

\bibitem{Artioli}
E.~Artioli, L.~Beir\~{a}o~da Veiga, and M.~Verani.
\newblock An adaptive curved virtual element method for the statistical
  homogenization of random fibre-reinforced composites.
\newblock {\em Finite Elem. Anal. Des.}, 177:103418, 12, 2020.

\bibitem{Voronoi2}
F.~Aurenhammer.
\newblock Voronoi diagrams—a survey of a fundamental geometric data
  structure.
\newblock {\em ACM Computing Surveys (CSUR)}, 23(3):345--405, 1991.

\bibitem{baines1994moving}
M.J. Baines.
\newblock {\em Moving finite elements}.
\newblock Oxford University Press, Inc., 1994.

\bibitem{MMFEM}
M.J. Baines, M.E. Hubbard, and P.K. Jimack.
\newblock A moving mesh finite element algorithm for the adaptive solution of
  time-dependent partial differential equations with moving boundaries.
\newblock {\em Applied Numerical Mathematics}, 54(3-4):450--469, 2005.

\bibitem{MovingReview}
M.J. Baines, M.E. Hubbard, and P.K. Jimack.
\newblock Velocity-based moving mesh methods for nonlinear partial differential
  equations.
\newblock {\em Communications in Computational Physics}, 10(3):509--576, 2011.

\bibitem{MMFEMscale}
M.J. Baines, M.E. Hubbard, P.K. Jimack, and A.C. Jones.
\newblock Scale-invariant moving finite elements for nonlinear partial
  differential equations in two dimensions.
\newblock {\em Applied Numerical Mathematics}, 56(2):230--252, 2006.

\bibitem{MovingPhaseChange}
M.J. Baines, M.E. Hubbard, P.K. Jimack, and R.~Mahmood.
\newblock A moving-mesh finite element method and its application to the
  numerical solution of phase-change problems.
\newblock {\em Communications in Computational Physics}, 6(3):595--624, 2009.

\bibitem{phasechange}
G.~Beckett, J.A. Mackenzie, and M.L. Robertson.
\newblock A moving mesh finite element method for the solution of
  two-dimensional stefan problems.
\newblock {\em Journal of Computational Physics}, 168(2):500--518, 2001.

\bibitem{Polycurved}
L.~Beir\~{a}o~da Veiga, F.~Brezzi, L.D. Marini, and A.~Russo.
\newblock Polynomial preserving virtual elements with curved edges.
\newblock {\em Math. Models Methods Appl. Sci.}, 30(8):1555--1590, 2020.

\bibitem{basicprinciples}
L.~Beir{\~a}o~da Veiga, F.~Brezzi, A.~Cangiani, G.~Manzini, L.D. Marini, and
  A.~Russo.
\newblock Basic principles of virtual element methods.
\newblock {\em Mathematical Models and Methods in Applied Sciences},
  23(01):199--214, 2013.

\bibitem{hitchiker}
L.~Beir{\~a}o~da Veiga, F.~Brezzi, L.D. Marini, and A.~Russo.
\newblock The hitchhiker's guide to the virtual element method.
\newblock {\em Mathematical models and methods in applied sciences},
  24(08):1541--1573, 2014.

\bibitem{ellipticVEM}
L.~Beir{\~a}o~da Veiga, F.~Brezzi, L.D. Marini, and A.~Russo.
\newblock Virtual element method for general second-order elliptic problems on
  polygonal meshes.
\newblock {\em Mathematical Models and Methods in Applied Sciences},
  26(04):729--750, 2016.

\bibitem{Bonito2013apriori}
A.~Bonito, I.~Kyza, and R.H. Nochetto.
\newblock {Time-discrete higher order ALE formulations: A priori error
  analysis}.
\newblock {\em Numerische Mathematik}, 125(2):225--257, 2013.

\bibitem{Bonito2013}
A.~Bonito, I.~Kyza, and R.H. Nochetto.
\newblock {Time-discrete higher-order ale formulations: Stability}.
\newblock {\em SIAM Journal on Numerical Analysis}, 51(1):577--604, 2013.

\bibitem{brenner2018virtual}
S.C. Brenner and L.~Sung.
\newblock Virtual element methods on meshes with small edges or faces.
\newblock {\em Mathematical Models and Methods in Applied Sciences},
  28(07):1291--1336, 2018.

\bibitem{Brezzi_DG_Darcy}
F.~Brezzi, T.J. Hughes, L.D. Marini, and A.~Masud.
\newblock Mixed discontinuous galerkin methods for darcy flow.
\newblock {\em Journal of Scientific Computing}, 22(1):119--145, 2005.

\bibitem{MFDcurved06}
F.~Brezzi, K.~Lipnikov, and M.~Shashkov.
\newblock Convergence of mimetic finite difference method for diffusion
  problems on polyhedral meshes with curved faces.
\newblock {\em Math. Models Methods Appl. Sci.}, 16(2):275--297, 2006.

\bibitem{blowup}
C.J. Budd, W.~Huang, and R.D. Russell.
\newblock Moving mesh methods for problems with blow-up.
\newblock {\em SIAM Journal on Scientific Computing}, 17(2):305--327, 1996.

\bibitem{NonHeirVEM}
A.~Cangiani, E.H. Georgoulis, and O.J. Sutton.
\newblock Adaptive non-hierarchical {G}alerkin methods for parabolic problems
  with application to moving mesh and virtual element methods.
\newblock {\em Math. Models Methods Appl. Sci.}, 31(4):711--751, 2021.

\bibitem{cangiani2017conforming}
A.~Cangiani, G.~Manzini, and O.J. Sutton.
\newblock Conforming and nonconforming virtual element methods for elliptic
  problems.
\newblock {\em IMA Journal of Numerical Analysis}, 37(3):1317--1354, 2017.

\bibitem{GCL}
W.~Cao, W.~Huang, and R.D. Russell.
\newblock A moving mesh method based on the geometric conservation law.
\newblock {\em SIAM Journal on Scientific Computing}, 24(1):118--142, 2002.

\bibitem{Bend3d}
F.~Dassi, A.~Fumagalli, A.~Scotti, and G.~Vacca.
\newblock Bend 3d mixed virtual element method for darcy problems.
\newblock {\em Computers \& Mathematics with Applications}, 119:1--12, 2022.

\bibitem{ALEreview}
J.~Donea, A.~Huerta, J.~Ponthot, and A.~Rodr{\'\i}guez-Ferran.
\newblock Arbitrary lagrangian-eulerian methods.
\newblock {\em Encyclopedia of Computational Mechanics Second Edition}, pages
  1--23, 2017.

\bibitem{lloydconvergence}
Q.~Du, M.~Emelianenko, and L.~Ju.
\newblock Convergence of the {L}loyd algorithm for computing centroidal
  {V}oronoi tessellations.
\newblock {\em SIAM journal on numerical analysis}, 44(1):102--119, 2006.

\bibitem{Ferguson22}
J.A. Ferguson, J.~Kópházi, and M.D. Eaton.
\newblock Nurbs enhanced virtual element methods for the spatial discretization
  of the multigroup neutron diffusion equation on curvilinear polygonal meshes.
\newblock {\em Journal of Computational and Theoretical Transport},
  51(4):145--204, 2022.

\bibitem{ALEmovingVoronoi}
E.~Gaburro, W.~Boscheri, S.~Chiocchetti, C.~Klingenberg, V.~Springel, and
  M.~Dumbser.
\newblock High order direct arbitrary-lagrangian-eulerian schemes on moving
  voronoi meshes with topology changes.
\newblock {\em Journal of Computational Physics}, 407:109167, 2020.

\bibitem{gelinas1981moving}
R.J. Gelinas, S.K. Doss, and K.~Miller.
\newblock The moving finite element method: applications to general partial
  differential equations with multiple large gradients.
\newblock {\em Journal of Computational Physics}, 40(1):202--249, 1981.

\bibitem{hallquist1985sliding}
J.O. Hallquist, G.L. Goudreau, and D.J. Benson.
\newblock Sliding interfaces with contact-impact in large-scale lagrangian
  computations.
\newblock {\em Computer methods in applied mechanics and engineering},
  51(1-3):107--137, 1985.

\bibitem{MovingMeshesBook}
W.~Huang and R.D. Russell.
\newblock {\em Adaptive moving mesh methods}, volume 174.
\newblock Springer Science \& Business Media, 2010.

\bibitem{MMPDEAnisotropic}
W.~Huang and Y.~Wang.
\newblock Anisotropic mesh quality measures and adaptation for polygonal
  meshes.
\newblock {\em Journal of Computational Physics}, page 109368, 2020.

\bibitem{MMFEMBC}
M.E. Hubbard, M.J. Baines, and P.K. Jimack.
\newblock Consistent dirichlet boundary conditions for numerical solution of
  moving boundary problems.
\newblock {\em Applied Numerical Mathematics}, 59(6):1337--1353, 2009.

\bibitem{jimack&Wathen}
P.K. Jimack and A.J. Wathen.
\newblock Temporal derivatives in the finite-element method on continuously
  deforming grids.
\newblock {\em SIAM Journal on Numerical Analysis}, 28(4):990--1003, 1991.

\bibitem{dGALEmapping}
K.~Lipnikov and N.~Morgan.
\newblock A high-order conservative remap for discontinuous galerkin schemes on
  curvilinear polygonal meshes.
\newblock {\em Journal of Computational Physics}, 399:108931, 2019.

\bibitem{MMFEMmonitors}
R.~Marlow, M.E. Hubbard, and P.K. Jimack.
\newblock Moving mesh methods for solving parabolic partial differential
  equations.
\newblock {\em Computers \& Fluids}, 46(1):353--361, 2011.

\bibitem{Masud_Darcy}
A.~Masud and T.J. Hughes.
\newblock A stabilized mixed finite element method for darcy flow.
\newblock {\em Computer methods in applied mechanics and engineering},
  191(39-40):4341--4370, 2002.

\bibitem{SedimentaryFlowVEM}
A.~Mazzia, M.~Ferronato, P.~Teatini, and C.~Zoccarato.
\newblock Virtual element method for the numerical simulation of long-term
  dynamics of transitional environments.
\newblock {\em Journal of Computational Physics}, page 109235, 2020.

\bibitem{miller1981moving}
K.~Miller and R.N. Miller.
\newblock Moving finite elements. i.
\newblock {\em SIAM Journal on Numerical Analysis}, 18(6):1019--1032, 1981.

\bibitem{MathBioBook}
J.D. Murray.
\newblock {\em Mathematical Biology: 1. An Introduction (3rd edition)}.
\newblock Springer-Verlag New York, 2002.

\bibitem{MovingHessian}
C.~Ngo and W.~Huang.
\newblock A study on moving mesh finite element solution of the porous medium
  equation.
\newblock {\em Journal of Computational Physics}, 331:357--380, 2017.

\bibitem{Voronoi1}
A.~Okabe, B.~Boots, K.~Sugihara, and S.N. Chiu.
\newblock {\em Spatial tessellations: concepts and applications of {V}oronoi
  diagrams}, volume 501.
\newblock John Wiley \& Sons, 2009.

\bibitem{Richter2017}
T.~Richter.
\newblock {\em {Fluid-structure Interactions}}, volume 118 of {\em Lecture
  Notes in Computational Science and Engineering}.
\newblock Springer International Publishing, Cham, 2017.

\bibitem{VEMMATLAB}
O.J. Sutton.
\newblock The virtual element method in 50 lines of matlab.
\newblock {\em Numerical Algorithms}, 75(4):1141--1159, 2017.

\bibitem{polymesher}
C.~Talischi, G.H. Paulino, A.~Pereira, and I.F.M. Menezes.
\newblock Polymesher: a general-purpose mesh generator for polygonal elements
  written in matlab.
\newblock {\em Structural and Multidisciplinary Optimization}, 45(3):309--328,
  2012.

\bibitem{hyperVEM}
G.~Vacca.
\newblock Virtual element methods for hyperbolic problems on polygonal meshes.
\newblock {\em Computers \& Mathematics with Applications}, 74(5):882--898,
  2017.

\bibitem{Vacca2018}
G.~Vacca.
\newblock An h-1-conforming virtual element for darcy and brinkman equations.
\newblock {\em Mathematical Models \& Methods in Applied Sciences},
  28(1):159--194, Jan 2018.

\bibitem{parabolicVEM}
G.~Vacca and L.~Beir{\~a}o~da Veiga.
\newblock Virtual element methods for parabolic problems on polygonal meshes.
\newblock {\em Numerical Methods for Partial Differential Equations},
  31(6):2110--2134, 2015.

\bibitem{PME}
J.L. V{\'a}zquez.
\newblock {\em The porous medium equation: mathematical theory}.
\newblock Oxford University Press, 2007.

\bibitem{Wang2019}
G.~Wang, F.~Wang, L.~Chen, and Y.~He.
\newblock A divergence free weak virtual element method for the stokes-darcy
  problem on general meshes.
\newblock {\em Computer Methods in Applied Mechanics and Engineering},
  344:998--1020, Feb 1 2019.

\bibitem{CFD}
P.~Wesseling.
\newblock {\em Principles of computational fluid dynamics}, volume~29.
\newblock Springer Science \& Business Media, 2009.

\bibitem{wriggers2006computational}
P.~Wriggers and T.A. Laursen.
\newblock {\em Computational contact mechanics}, volume~2.
\newblock Springer, 2006.

\bibitem{VEMcontact}
P.~Wriggers, W.T. Rust, and B.D. Reddy.
\newblock A virtual element method for contact.
\newblock {\em Computational Mechanics}, 58(6):1039--1050, 2016.

\bibitem{Zhao2020}
J.~Zhao, B.~Zhang, S.~Mao, and S.~Chen.
\newblock The nonconforming virtual element method for the darcy-stokes
  problem.
\newblock {\em Computer Methods in Applied Mechanics and Engineering}, 370, Oct
  1 2020.

\end{thebibliography}
\end{document}